\documentclass{article}

\usepackage{arxiv}

\usepackage[utf8]{inputenc} 
\usepackage[T1]{fontenc}    
\usepackage{hyperref}       
\usepackage{url}            
\usepackage{booktabs}       
\usepackage{amsfonts}       
\usepackage{nicefrac}       
\usepackage{microtype}      
\usepackage{lipsum}		
\usepackage{amsthm}
\usepackage{graphicx}
\usepackage{doi}
\usepackage{physics}
\usepackage{amssymb}
\usepackage{pifont}
\usepackage{algorithm}
\usepackage{algpseudocode}
\usepackage{titlesec} 
\usepackage{multirow}
\theoremstyle{remark}
\newtheorem{remark}{Remark}

\newcommand{\pd}[2]{\frac{\partial #1}{\partial #2}}
\newcommand{\pdd}[2]{\frac{\partial^2 #1}{\partial #2^2}}
\newcommand{\Ab}{\mathbf{A}}
\newcommand{\Abt}{\tilde{\mathbf{A}}}
\newcommand{\ab}{\mathbf{a}}
\newcommand{\abt}{\tilde{\mathbf{a}}}
\newcommand{\fb}{\mathbf{f}}
\newcommand{\fsindy}{\hat{\mathbf{f}}}

\newcommand{\mb}{\mathbf{m}}
\newcommand{\vbb}{\mathbf{v}}
\newcommand{\Gb}{\mathbf{G}}

\usepackage{setspace}
\title{Learning gradient flow: using equation discovery to accelerate engineering optimization}

\date{} 					

\author{
    Grant Norman\thanks{These authors contributed equally.} \\
	Smead Aerospace Engineering Sciences\\
    3775 Discovery Drive \\
	Boulder, CO 80309 \\
	\texttt{grant.norman@colorado.edu} \\
    \And
    Conor Rowan\footnotemark[1] \\
	Smead Aerospace Engineering Sciences\\
    3775 Discovery Drive \\
	Boulder, CO 80309 \\
	\texttt{conor.rowan@colorado.edu} \\
    \And
    Kurt Maute \\
    Smead Aerospace Engineering Sciences\\
    3775 Discovery Drive \\
	Boulder, CO 80309 \\  
    \texttt{kurt.maute@colorado.edu} \\
    \And
    Alireza Doostan \\
    Smead Aerospace Engineering Sciences\\
    3775 Discovery Drive \\
	Boulder, CO 80309 \\
	\texttt{alireza.doostan@colorado.edu} \\
}

\renewcommand{\headeright}{}
\renewcommand{\undertitle}{}
\renewcommand{\shorttitle}{}

\hypersetup{
pdftitle={Learning gradient flow: using equation discovery to accelerate engineering optimization},
pdfsubject={cs.NA, math.DS, math.OC},
pdfauthor={Grant~Norman, Conor~Rowan, Alireza~Doostan, Kurt~Maute},
pdfkeywords={Equation discovery, data-driven modeling, gradient flow, surrogate models, optimization},
}

\begin{document}
\maketitle

\vspace{-0.5cm}
\begin{abstract}

In this work, we investigate the use of data-driven equation discovery for dynamical systems to model and forecast continuous-time dynamics of unconstrained optimization problems. To avoid expensive evaluations of the objective function and its gradient, we leverage trajectory data on the optimization variables to learn the continuous-time dynamics associated with gradient descent, Newton's method, and ADAM optimization. The discovered gradient flows are then solved as a surrogate for the original optimization problem. To this end, we introduce the Learned Gradient Flow (LGF) optimizer, which is equipped to build surrogate models of variable polynomial order in full- or reduced-dimensional spaces at user-defined intervals in the optimization process. We demonstrate the efficacy of this approach on several standard problems from engineering mechanics and scientific machine learning, including two inverse problems, structural topology optimization, and two forward solves with different discretizations. Our results suggest that the learned gradient flows can significantly expedite convergence by capturing critical features of the optimization trajectory while avoiding expensive evaluations of the objective and its gradient.

\end{abstract}



\newpage

\section{Introduction}
\label{sec:intro}

Optimization is a pillar of engineering design and analysis, as well as computational science writ large. As traditionally presented, optimization problems entail finding a set of inputs, or optimization variables, to a scalar objective function which corresponds to a minimum value of that function \cite{ruszczynski_nonlinear_2006}.
This simple, yet powerful idea finds use in diverse applications across disciplines: partial differential equation (PDE) solutions, eigenvalue problems, optimal control, engineering design, inverse problems, neural network training, and many other tasks in scientific computing; see, e.g., \cite{e_deep_2017, rowan_solving_2025, brogan_modern_1991, zhang_non-linear_2024, raissi_physics-informed_2019, young_fully_2024}.
Because most optimization problems do not have analytical solutions, a common numerical approach is to start with an initial guess and to iteratively evolve the optimization variables in a direction that is expected to decrease the value of the objective. Different optimization algorithms use different techniques to estimate a direction of decreasing objective value. With gradient descent, the gradient of the objective function at the current state of the optimization variables defines the update direction. With a Newton-type optimization scheme, a quadratic approximation of the objective is built around the current iterate, and the update is taken in the direction of the minimum of this quadratic approximation \cite{ruszczynski_nonlinear_2006}. In the machine learning community, there are many variants of gradient descent, such as the now-famous ``ADAM'' algorithm \cite{kingma_adam_2014}, which leverages first- and second-moment estimates of past gradients to scale the update direction and improve convergence. All of these algorithms have a notion of an optimal search direction and a ``step size'' or ``learning rate,'' which controls the size of the update to the optimization variables at each iteration. Note that a systematic comparison of some popular second-order optimizers was recently conducted in \cite{kiyani_optimizing_2025}.
While many optimization strategies do not require gradient information \cite{larson_derivative-free_2019}, we restrict our attention to gradient-based optimization.

In the limit as the step size approaches zero, a continuous-time system of ordinary differential equations for the optimization variables is obtained. For example, gradient descent updates resemble a forward Euler strategy for time integrating the optimization variables. Thus, in the continuous time limit of gradient descent, the optimization problem corresponds to ``gradient flow'' of the objective function \cite{lemarechal_cauchy_2012}. Stable attractors of this dynamical system are solutions to the optimization problem, and the standard discrete-time gradient descent scheme is returned when the time derivative is approximated with finite differencing. Noting this equivalence between optimization and dynamical systems, many authors have used convergence, stability, and control analyses from the study of ordinary differential equations (ODEs) in order to understand the properties of optimization algorithms. The neural tangent kernel (NTK)---a continuous time formulation of neural network training that linearizes the network parameters around the origin---has been studied extensively to understand the dynamics of training \cite{jacot_neural_2018}. In \cite{kumar_towards_2025}, the authors extend the NTK analysis of training dynamics to parameters that are far from the origin. With the help of continuous time gradient flows, convergence guarantees for training neural networks with rectified linear unit (ReLU) activation functions were given in \cite{jentzen_convergence_2023}.
An analysis of the stochastic differential equation corresponding to mini-batched gradient descent is performed in \cite{orvieto_continuous-time_2020, latz_analysis_2021} to show convergence properties for certain non-convex objective functions.
Various stability analyses of neural network training have been conducted using the gradient flow dynamics \cite{rosca_continuous_2023, garg_fixed-time_2021}. Drawing from control theory, the authors in \cite{chen_accelerated_2024}  use a proportional-integral-derivative controller on the gradient flow dynamics in order to more aggressively explore the objective landscape.
The analysis of continuous time gradient flow is extended to the case of equality-constrained objectives in \cite{chen_sign_2020}. Furthermore, the authors in \cite{kovachki_continuous_2021} show that optimization algorithms with history dependence (e.g., momentum) can also be understood in continuous time. Furthermore, the Nesterov accelerated gradient method can be shown to have a continuous time analogue, which provides insight into its behavior by establishing a clear connection to physical dynamical systems \cite{su_differential_2015}. For a comprehensive perspective on the connection between accelerated and stochastic optimization and dissipative systems, the reader is referred to \cite{franca_gradient_2021}. By way of example, the continuous time Nesterov method has been shown to improve the training of autoencoders \cite{flouris_explicit_2025}. Similarly, continuous time gradient flows are used in the volume reconstruction literature in order to find minima of objective functionals penalizing the curvature or surface area of the reconstructed volume \cite{bretin_learning_2022, bretin_penalized_2024}. 

The continuous time perspective offers a useful connection between dynamical systems and optimization, but does not offer a remedy for the inherent cost of gradient-based optimization.
All such methods require computing at least the objective and its gradient for every iteration, and they may take hundreds or thousands of iterations to converge. In the context of engineering optimization, this often necessitates repeated solves of large-scale PDE systems. Thus, owing to this high computational cost, there is significant interest in devising techniques to accelerate optimization problems. Several strategies have been proposed to this end. Surrogate modeling is a popular technique to replace expensive aspects of the optimization process with cheap-to-evaluate mathematical functions that are fit from data. Some surrogate models replace the entire optimization process, going directly from inputs to the optimization problem to the final solution. For example, in \cite{anantha_padmanabha_solving_2021}, a whole class of inverse problems can be solved efficiently by a surrogate that learns the mapping between an output field and the corresponding material property. Similarly, works such as \cite{barmada_deep_2021, giacomini_surrogate_2025} aim to replace an entire structural topology optimization pipeline with a trained surrogate that maps problem parameters to the final design. Other surrogate models approximate the dependence of the loss on the optimization variables and use this approximation to compute gradients. Though Gaussian process surrogates (sometimes called Kriging models) are usually gradient-free, in \cite{ozkaya_gradient_2024}, a Kriging surrogate for the objective function is augmented with gradient information around observation points in the context of aerodynamic shape optimization. In \cite{backhaus_gradient_2013}, neural networks are used to approximate the dependence of the loss on the design parameters in order to avoid expensive online forward solves during optimization. We note that even standard second-order optimization methods, such as Newton and quasi-Newton methods, can be conceptualized as constructing local surrogate models of the objective in order to expedite convergence. In a standard Newton approach to finding the stationary point of an objective function, a quadratic surrogate of the objective is made around the current state of the design variables and is used to determine optimal step sizes and directions \cite{ruszczynski_nonlinear_2006}. Similarly, with quasi-Newton methods such as BFGS, data from the optimization history is used on an ongoing basis to estimate this quadratic approximation \cite{nocedal_numerical_2006}. Though not typically conceptualized in this way, we suggest that even quasi-Newton methods somewhat resemble data-driven surrogates.

When the surrogate model approximates a mapping from one function space to another, ``operator learning'' frameworks such as the Deep Operator Network (DeepONet, \cite{lu_deeponet_2021}), Fourier Neural Operators (FNO, \cite{li_fourier_2021}), or the Laplace Neural Operator (LNO, \cite{cao_lno_2023}) are often employed. For example, pre-trained surrogate models of this type can be used to make informed initializations for optimization-based PDE solvers, which dramatically decrease the number of iterations to convergence \cite{aghili_accelerating_2025, abdelrehim_active_2025, fesefeldt_using_2024}.
An alternative to reducing the number of steps to convergence is to reduce the cost per iteration.
On this front, surrogate models have been used to approximate gradients of the loss function based on data seen in an offline training stage.
A non-intrusive approach to approximating gradients in the context of inverse problems is proposed in \cite{xiao_deep-learning-based_2021}. Other examples of data-driven surrogates replacing expensive gradient computations are given in \cite{ozan_data-driven_2024, xu_machine_2020}. Some surrogates for objective gradients have been designed to learn gradient information online, reducing the overhead cost of offline training \cite{xing_accelerating_2025}.

While pre-trained surrogates can reduce the number of optimization steps with informed initializations, another strategy to improve performance---both in the sense of the number of steps to converge and the quality of the converged solution---is to devise better optimization algorithms. For example, the popular ADAM algorithm has shown that a carefully constructed update strategy can lead to improved performance even on highly non-convex and stochastic objectives \cite{kingma_adam_2014}. Historically, optimization algorithms have been designed and tuned manually. Recognizing this, a number of authors have explored data-driven automation of the design of optimization algorithms. In the ``Learning to Optimize'' framework, the iterative update to the optimization variables is represented as a neural network with inputs given by the objective function and its derivatives \cite{li_learning_2016}. Reinforcement learning or gradient-based meta-learning is used to learn the optimal update strategy based on optimization dynamics observed from a given class of problems. This methodology was extended and validated in a series of follow-up works \cite{chen_learning_2021, yin_learning_2022}. A similar strategy is proposed in \cite{andrychowicz_learning_2016}, where the optimization update rule is parameterized and learned for a class of problems.

A parallel line of work, ostensibly unrelated to optimization, is the discovery of dynamical systems from data using symbolic regression and equation discovery. Early work on the discovery of laws from data using symbolic regression dates back two decades \cite{bongard_automated_2007}, but gained traction in 2016 with the Sparse Identification of Nonlinear Dynamics (SINDy) framework for equation discovery \cite{brunton_discovering_2016}. SINDy takes the measured trajectory of a system governed by a nonlinear system of ordinary differential equations and returns a sparse approximation of its dynamics. We note that equation discovery uses sparse regression techniques to find an optimal approximation of the system dynamics, whereas symbolic regression uses genetic algorithms to search through a large space of symbolic expressions. An alternative equation discovery framework to the SINDy methodology is ``Neural Ordinary Differential Equations'' (Neural ODEs), which represents the unknown dynamics as a neural network taking in the system state and its derivatives \cite{chen_neural_2019}. This method has also been studied in a variety of contexts, such as model reduction \cite{caldana_neural_2024}, fluid mechanics \cite{wang_learning_2024}, and reactor systems \cite{sorourifar_physics-enhanced_2023}. We remark that Gaussian process regression has also been used for the discovery of governing equations by providing uncertainty-aware estimates of the system state and its derivatives, which can then be used to identify terms in a differential equation \cite{raissi_machine_2017}. In this work, our focus will be on equation discovery with SINDy, though the methods reviewed above could be used as well. The SINDy methodology has been successfully applied to problems in biology \cite{massonis_distilling_2023}, to systems with noisy measurements \cite{messenger_weak_2021,hokanson_simultaneous_2022,wentz_derivative-based_2023-1}, and to systems whose dynamics are describable in terms of low-dimensional subspaces \cite{champion_data-driven_2019}.

Once the connection between optimization and dynamical systems is established, obtaining a surrogate model for the gradients of the objective becomes equivalent to learning a dynamical system. In this work, we aim to bring together the dynamical perspective on optimization, the need for surrogate models to accelerate optimization problems, and ideas from data-driven non-linear equation discovery, here SINDy. Namely, we apply a SINDy-inspired approach while training in order to learn the dynamics of optimization based on past observations of the loss function and its gradient. This SINDy model, built from a user-defined subset of the optimization history, can then be used as a surrogate for the full loss and gradient evaluations to forecast a user-defined number of steps. Our contributions can be broken down as follows:

\begin{enumerate}
    \item We propose the use of SINDy-based surrogates of the continuous-time dynamics of optimization,
    \item We extend this methodology to build surrogates for gradient descent, Newton, and ADAM optimization,
    \item We develop a PyTorch optimizer that builds SINDy models from the optimization trajectory data in real time,
    \item We showcase the ability of our method to accelerate optimization on five standard problems from engineering optimization and scientific machine learning.
\end{enumerate}

The rest of this paper is organized as follows. In Section 2, we introduce the SINDy methodology and discuss how it can be adapted to be used online as a surrogate for the optimization dynamics. In Section 3, we describe our custom optimization framework called Learned Gradient Flow (LGF). In Section 4, we present five numerical examples designed to cover a range of problems from engineering mechanics, including i) linear and nonlinear physics, ii) spectral and finite element discretizations, and iii) forward solves, inverse problems, and topology optimization. In Section 5, we close with concluding remarks and propose directions for future work.


\section{Background}


%
%
\subsection{Modeling ODEs through SINDy}
\label{sec:sindy_for_ode_discovery}

Here, we introduce notation and basic details for SINDy \cite{brunton_discovering_2016}.
As previously mentioned, SINDy methods construct models for ODE systems based on observations of the state, often involving a relatively cheap solve of a convex program.
To be more concrete, consider a state vector $\ab(t) \in \mathbb{R}^n$ which evolves according to some (potentially unknown) ordinary differential equation
$$
\frac{\dd \ab}{ \dd t} = \fb ( \ab ) .
$$
Assume we measure this state $k$ times, yielding a measurement matrix
$$
\Ab = 
\begin{bmatrix}
    \ab (t_1) & \ab (t_2) & \dots & \ab (t_k) \\
\end{bmatrix}^T \in \mathbb{R}^{k \times n}.
$$
Through either numerical differentiation or access to $\fb$ itself, assume we also have the time derivative of the measurement matrix,
$$
\dot{\Ab} = 
\begin{bmatrix}
    \frac{\dd \ab}{ \dd t} (t_1)  & \frac{\dd \ab}{ \dd t} (t_2) & \dots & \frac{\dd \ab}{ \dd t} (t_k) \\
\end{bmatrix} ^T \in \mathbb{R}^{k \times n}.
$$
Consider a candidate library function $\boldsymbol{\phi} : \mathbb{R}^n \to \mathbb{R}^{p}$.
Each element of this function corresponds to one of the $p$ candidate functions in a pre-defined library.
For instance, denoting the elements of $\ab$ as $a_1, a_2, \dots, a_n$, a candidate function library including up to quadratic polynomial terms would be
$$
\boldsymbol{\phi} (\ab)^T = \begin{bmatrix}
    1 &
    a_1 &
    a_2 &
    \dots &
    a_n &
    a_{1}^{2} & 
    a_{1} a_{2} &
    \dots &
    a_{n}^{2}
\end{bmatrix} .
$$
In some instances, prior knowledge is used to select a specific candidate function library (e.g., trigonometry functions), but in this work, we assume only polynomial libraries, including cross terms, up to some total degree $\mathcal{P}$.
The total number of library elements can then be computed as
\begin{equation}
    p =\binom{n + \mathcal{P}}{\mathcal{P}} .
    \label{eq:polynomial_size_scaling}
\end{equation}

Repeatedly applying these (polynomial) candidate functions at different time steps gives the library \textit{matrix}:
$$
\mathbf{\Phi} (\Ab) = 
\begin{bmatrix}
    \boldsymbol{\phi} (\ab (t_1))^T \\
    \boldsymbol{\phi} (\ab (t_2))^T \\
    \vdots \\
    \boldsymbol{\phi} (\ab (t_k))^T \\
\end{bmatrix} \in \mathbb{R}^{k \times p},
$$
SINDy seeks a coefficient matrix $\mathbf{\Xi} \in \mathbb{R}^{p \times n}$ such that
\begin{equation}
    \dot{\Ab} \approx \mathbf{\Phi} (\Ab) \mathbf{\Xi},
    \label{eq:sindy_strong_form_goal}
\end{equation}
and $\mathbf{\Xi}$ is sparse.
In the SINDy literature, much emphasis is placed on this sparsity, since the goal is often to ``interpret'' the discovered ODE, giving rise to a multitude of methods such as sequentially thresholded least squares, stepwise sparse regression, sparse relaxed regularized regression, and forward regression orthogonal least squares \cite{kaptanoglu_pysindy_2022}
\footnote{\cite{rowan_definition_2025} argues that such sparsity is insufficient for interpretability, and that interpretability should instead emphasize mechanistic understanding.}.
In the present work, on the other hand, we are not concerned with interpretability, as we will only use the discovered ODE to generate predictions.
Thus, we opt for a least squares solution, but note that sparsity-promoting methods may act as a means to regularize the problem of determining $\mathbf{\Xi}$ \cite{candes_stable_2005}, constituting a possible future research direction.

Other SINDy works prioritize a weak form of Eq. \eqref{eq:sindy_strong_form_goal}, shifting the time derivative from the potentially noisy state variables to analytic test functions by applying integration by parts \cite{reinbold_using_2020, messenger_weak_2021-1}.
We refer to this class of Galerkin variants broadly as ``WSINDy.''
Once $\mathbf{\Xi}$ is computed by some combination of these aforementioned methods, new predictions can be generated by solving the ODE
$$
\frac{\dd \ab}{\dd t} (t) = \mathbf{\Xi}^T \boldsymbol{\phi} (\ab (t))  =: \fsindy(\ab(t)).
$$


\subsection{Other methods for modeling ODEs}
While the methods presented in Sec. \ref{sec:methods} will most closely rely upon SINDy, we also overview some related work and possible alternatives for learning ODEs from data, for the sake of completeness.
Thus, the task is still to learn an approximation $ \fsindy$ of $\fb$ such that 
$$
\frac{\dd \ab}{\dd t} = \fsindy(\ab(t)) ,
$$
based on the samples of the state variables, $\Ab$.
Rather than prescribing $\fsindy$ as a linear expansion of library functions, we could specify $\fsindy$ as a neural network, yielding a neural ODE \cite{chen_neural_2019}.
As a result of this nonlinear approximation, the parameters of $\fsindy$ (previously $\mathbf{\Xi}$) must now be determined through an iterative minimization process.
The objective is given as
$$
\sum_{i=1}^{k} \norm{
        \ab (t_i) - \int_{0}^{t_i} \fsindy ( \ab ( \tau ) ) \, \dd \tau 
    }_{2}^{2} .
$$
While avoiding numerical differentiation of potentially noisy data, these methods require numerically integrating the state over time. Further, the nonlinear form of $\fsindy$ requires an iterative solver, which is often more expensive than finding the coefficients in a linear representation through a direct solve (if we are not concerned with sparsity). Other autoregressive models or flow maps instead directly learn the mapping from the state at one time to the next, thus omitting an explicit numerical integration scheme \cite{churchill_flow_2023, brandstetter_message_2023}.
Denoting such a model as $\hat{\mathbf{g}} (\ab)$ and repeated applications with an exponent, the equivalent objective function is
%
$$
\sum_{i=1}^{k} \norm{
    \ab (t_i)
    - \hat{\mathbf{g}}^{i} (\ab (0) )
} .
$$
%
While the contributions presented in Sec. \ref{sec:learned_gradient_flow} could be modified to use either Neural ODEs, autoregressive models, or Gaussian processes, we use a technique similar to SINDy due to its reduced computational cost and clearer connection to ODEs.

\subsection{Gradient flow ODE}
\label{sec:gradient_flow_ode}

Here, we review an idea mentioned in Sec. \ref{sec:intro}, namely that an ODE system emerges when the learning rate of gradient descent is taken to zero. For a learning rate $\eta > 0$, an objective function $z: \mathbb{R}^n \to \mathbb{R}$, and current optimization variables $\ab_k$ at iteration $k$, gradient descent computes the optimization variables at iteration $k+1$ as 
\begin{equation}
    \ab_{k+1} = \ab_k - \eta \cdot \pd{}{\mathbf{a}} \, z (\ab_k).
    \label{eq:gradient_descent}
\end{equation}
By noting the similarity to the forward Euler derivative approximation and taking $\eta \to 0$, we have
\begin{equation}
    \frac{\dd \ab}{\dd t} (t_k) = -\pd{}{\mathbf{a}} \, z (\ab (t_k) ),
    \label{eq:gradient_flow}
\end{equation}
where $\ab (t_k) = \ab_k$.
Thus, rather than updating $\ab_k$ through discrete applications of Eq. \eqref{eq:gradient_descent}, we can introduce a continuous analogue and evolve $\ab (t)$ according to the ODE in Eq.~\eqref{eq:gradient_flow}.
The optimization variables throughout optimization iterations correspond to a solution of an ODE over time, where the iteration number and time value are related through the learning rate $\eta$:
$t_k = k \cdot \eta$.
Motivated by this, we sometimes refer to the optimization variables as the ``state variables,'' when considering this dynamical systems interpretation
\footnote{While the ODE presented in Eq.~\eqref{eq:gradient_flow} is first-order and autonomous, other optimization methods can yield different types of ODEs, such as a second-order, non-autonomous system in the case of Nesterov acceleration \cite{su_differential_2016}. The methods presented in Sec. \ref{sec:methods} can be adapted to other ODE forms, without loss of generalization.}. We remark that this continuous-time formulation of gradient descent does not allow for constraints on the objective function. Accordingly, we restrict our attention to unconstrained optimization problems through this work.



\section{Method}
\label{sec:methods}

In this section, we build on ideas from Sec. \ref{sec:sindy_for_ode_discovery} and \ref{sec:gradient_flow_ode} to give the core contribution of this work.

\subsection{Learned gradient flow}
\label{sec:learned_gradient_flow}

While we can treat an unconstrained optimization problem as a dynamical system, solving this dynamical system for the state is expensive due to the cost of evaluating $\partial z / \partial \mathbf{a}$.
Instead, we can record the optimization variables $\ab$ at various iterations during some standard training procedure (Eq. \eqref{eq:gradient_descent}) and use SINDy to build a data-driven model for the gradient flow equation, which we name $\fsindy$.
Then, rather than integrating the expensive, true gradient equation, we integrate this cheaper surrogate ODE to evolve $\ab$ over training iterations.
In other words, we evolve the optimization variables according to
\begin{equation}
    \frac{\dd \ab}{\dd t} (t) = \fsindy ( \ab (t) ),
    \label{eq:learned_gradient_flow}
\end{equation}
where we construct $\fsindy$ as
\begin{equation}
    \fsindy (\ab ) = \mathbf{\Xi}^T \boldsymbol{\phi} (\ab (t)),
\end{equation}
based on data generated by an initial training phase where we observe only instances of the optimization variables $\ab$.
That is, we do not record instances of the gradient, although this is an option we later adopt for ADAM in Sec.~\ref{sec:learned_adam}.

\begin{remark}
    \label{rem:rhs_grad_tracking}
    In the learned gradient flow framework, one could alternatively record gradient evaluations during the data-generating phase and fit $\fsindy$ directly to the right-hand side $-\pd{}{\mathbf{a}} z(\ab)$, in the same spirit as the gradient-modeling step we later adopt for ADAM in Sec.~\ref{sec:learned_adam}.
    For gradient descent with fixed step size $\eta$, this ``tracked-gradient'' approach is equivalent to inferring $\dot{\Ab}$ from just the trajectory, where we estimate $\dot{\Ab}$ via forward differences (see Eq. \eqref{eq:gradient_descent} and \eqref{eq:gradient_flow}).
    While we generally estimate $\dot{\Ab}$ using centered differences, in Sec.~\ref{sec:eresults:ex1_2var} we observed consistent behavior using forward, second-order centered, and fourth-order centered finite differences, possibly due to accurate underlying gradients (low noise).
    A systematic comparison of these alternatives (which include directly modeling the right-hand side) is beyond the scope of the present work.
\end{remark}



A successful construction of $\fsindy (\ab )$ leads to the surrogate $\fsindy (\ab )\approx -\pd{}{\mathbf{a}} \, z (\ab )$. 

Let $\mathcal K$ be the number of iterations of the standard training procedure for which the optimization variables were recorded.
Then, by specifying a total number of evolution iterations $\mathcal M > \mathcal K$ for the learned ODE, we also specify the final integration time $t = \mathcal M \cdot \eta$.
The larger the difference between $\mathcal M$ and $\mathcal K$, the greater the computational savings. We summarize these core steps in Algorithm \ref{alg:learned_gradient_flow}.
For the time integration, we use the \verb|torchdiffeq| library, generally with the default Dormand-Prince adaptive method (\verb|dopri5|) \cite{chen_torchdiffeq_2018}.

\begin{algorithm}
\caption{Learned Gradient Flow}\label{alg:learned_gradient_flow}
    \begin{algorithmic}[1]
        \Require $\ab_1$ \Comment{initial state}
        \Require $\mathcal K$ \Comment{number of true gradient descent iterations/state samples}
        \Require $\mathcal M \geq \mathcal K$ \Comment{number of total steps, including surrogate ``steps''}
        \Require $\eta$ \Comment{gradient descent learning rate, or time between state samples}
        \State $\Ab \gets \ab_1$
        \For{$k=1,\dots,\mathcal K-1$} \Comment{data generating phase}
        \State $\ab_{k+1}\gets \ab_k - \eta \cdot \partial  z(\ab_k) / \partial \mathbf{a}_k$ \Comment{standard gradient descent update}
        \State append $\ab_{k+1}$ as a column to $\Ab$
        \EndFor
        \State $\mathbf{\Phi} \gets \mathbf{\Phi}(\Ab)$
        \State Compute $\dot{\Ab}$ \Comment{or replace with WSINDy}
        \State Solve SINDy system for $\mathbf{\Xi}$
        \State Define $\fsindy (\ab)=\mathbf{\Xi}^T \boldsymbol{\phi} (\ab )$ \Comment{cheap to evaluate}
        \State Integrate $\frac{\dd \ab}{\dd t} = \fsindy (\ab)$ from $t = \mathcal K \cdot \eta$ to $t = \mathcal M \cdot \eta$ \Comment{use $\ab_{\mathcal  K}$ as the ``initial condition''}
        \State \Return $\ab (t_\mathcal M)$
    \end{algorithmic}
\end{algorithm}

In Sec. \ref{sec:methods:retraining}, we discuss these parameters $\mathcal  K$ and $\mathcal  M$ more and report them for the numerical experiments in Sec. \ref{sec:results}.
Generally, more complicated optimization dynamics (such as a larger $\ab$) require more training data to identify (hence a larger $\mathcal  K$).
We use only polynomial library functions for $\boldsymbol{\phi}$, and note that for our examples, a low total order suffices.
The size of $\boldsymbol{\phi} (\ab)$ scales poorly with $\ab$ and the total polynomial order $\mathcal{P}$ (see Eq. \eqref{eq:polynomial_size_scaling}), thus motivating low polynomial orders and/or applications with few design variables.
Alg. \ref{alg:learned_gradient_flow} is written in a general way to accommodate various SINDy methods.
In our examples, we use a simple variant where the derivatives of the state variables ($\dot{\Ab}$) are computed via second-order centered finite differences and the coefficients $\mathbf{\Xi}$ are determined through sequentially thresholded least squares applied to the fitting problem associated with Eq. \eqref{eq:sindy_strong_form_goal}.
Our solver resembles that of \texttt{STLSQ} in \texttt{PySINDy} where $\alpha$ specifies the $\ell_2$ regularization weight used in ridge regression, and coefficients less than some \texttt{threshold} are removed from the problem after each iteration of the solve.
In almost all of the examples we will present in Sec. \ref{sec:results}, we set $\alpha$ to $10^{-6}$ and \texttt{threshold} to $10^{-8}$, normalize columns of $\mathbf{\Phi} (\Ab)$, complete up to $20$ iterations (\texttt{max\_iter~=~20}), and perform a final unregularized solve to avoid returning a biased $\mathbf{\Xi}$ (\texttt{unbias = True}).
For these examples, we observe similar results when using vanilla least squares (i.e., $\alpha = 0$ and \texttt{threshold} = $0$).
This suggests that, for our applications, solver choices are secondary to the assumption that past optimization dynamics are representative of future dynamics over the prediction horizon.
That said, this insensitivity to $\alpha$ and \texttt{threshold} is not expected in general: ridge regularization and thresholding can substantially improve the solved $\mathbf{\Xi}$ when $\dot{\Ab}$ is noisy, the library is ill-conditioned, or candidate terms are strongly correlated.

Importantly, the key contribution of LGF is integrating a learned dynamical system within optimization, allowing for developments in data-driven ODE modeling to be repurposed for LGF without loss of generality.
For instance, in the case of stochastic or mini-batch gradient descent, the objective gradient is only estimated at each iteration, yielding a noisier optimization trajectory, and we hypothesize that in this case, the more noise-robust WSINDy or regularized (or sparsity-promoting) solution techniques could become increasingly beneficial.
More generally, this model-agnostic perspective emphasizes that our contribution is not tied to a particular SINDy instantiation, but rather a general framework in which the model identification step can be replaced by any procedure that infers a surrogate ODE model from the optimization history.




%
\subsection{Generalization to other optimizers}  

While we initially introduced our approach by using gradient descent as an example, other optimization methods such as ADAM and Newton's method are commonly used in practice \cite{boyd_convex_2004, kingma_adam_2017}.
We overview the generalization to those methods here.

\subsubsection{Newton's method}
\label{sec:learned_newtons}

For a damping parameter $\eta$, Newton's method is given as
\begin{equation}
    \ab_{k+1} = \ab_{k} - \eta \cdot \Big(   \frac{\partial^2 }{\partial \mathbf{a} \partial \mathbf{a}} \, z (\ab_k)\Big)^{-1} \pd{}{\mathbf{a}} \, z(\ab_{k}),
    \label{eq:newton_update_rule}
\end{equation}
where $\partial^2 z /\partial \mathbf{a} \partial \mathbf{a}$ is the Hessian matrix of the objective $z$ \cite{boyd_convex_2004}. The Newton update rule in Eq. \eqref{eq:newton_update_rule} is similar to gradient descent in that it lends itself to a straightforward dynamical systems interpretation. The continuous-time Newton method is given by 
\begin{equation}\label{newtonflow}
\frac{\mathrm{d} \ab}{\mathrm{d}t} = - \qty( \frac{\partial^2 z}{\partial \mathbf{a} \partial \mathbf{a}})^{-1} \pd{z}{\mathbf{a}}, \quad \mathbf{a}(0)=\mathbf{a}_0,
\end{equation}
where $\mathbf{a}_0$ is the initialization of the optimization variables. As a result, we apply the same procedure as in Algorithm \ref{alg:learned_gradient_flow}, simply replacing the gradient descent update in the data-generating phase with Newton updates.
Subsequently, we still integrate the optimization variables in time according to Eq. \eqref{eq:learned_gradient_flow}, wherein $\fsindy$ seeks to approximate the time derivative of $\ab$.
%
It is worth highlighting that generating $\fsindy$ still relies only on data from the state variables $\ab$ rather than data directly involving the Hessian.

\subsubsection{ADAM}
\label{sec:learned_adam}

ADAM adds both momentum and normalization to gradient descent, relying on two history variables for each parameter \cite{kingma_adam_2015}.
The final update is given as
\begin{equation}
    \ab_{k+1} = \ab_{k} - \eta \cdot \frac{\hat{\mb}_{k+1}}{\sqrt{\hat{\vbb}_{k+1}} + \epsilon},
    \label{eq:discrete_adam:param_update}
\end{equation}
where the square root and division are element-wise, and the unbiased momentum parameters are, respectively, computed for $0 \leq \beta_1, \beta_2 < 1$ as
\begin{equation}
    \begin{aligned}
        \hat{\mb}_{k+1} &= \frac{\mb_{k+1}}{1 - \beta_{1}^{k+1}}, \\
        \hat{\vbb}_{k+1} &= \frac{\vbb_{k+1}}{1 - \beta_{2}^{k+1}}.
    \end{aligned}
    \label{eq:discrete_adam:biased_moments}
\end{equation}    
The biased history variables, $\mb$ and $\vbb$, begin as zero and are iterated as
\begin{equation}
    \begin{aligned}
        \mb_{k+1} &= \beta_1 \mb_{k} + (1 - \beta_1) \pd{}{\mathbf{a}} z (\ab_{k}) , \\
        \vbb_{k+1} &= \beta_2 \vbb_{k} + (1 - \beta_2) \left ( \pd{}{\mathbf{a}} z (\ab_{k}) \right )^2,
    \end{aligned}
    \label{eq:discrete_adam:unbiased_update}
\end{equation}
where the square is element-wise. Importantly, $\ab_{k+1}$ depends not only on $\ab_{k}$ as in gradient descent and Newton's method but also on $\ab_{k-1}, \dots, \ab_{0}$, through the history variables.
In modeling the evolution of $\ab (t)$ with a differential equation, we must therefore either use a delay differential equation or also include $\mb (t)$ and $\vbb (t)$ in the system.
We opt for the latter approach, yielding
\begin{equation}
    \begin{aligned}
        \frac{\dd \ab}{\dd t} (t) &= - \frac{\frac{\mb (t)}{1 - \beta_{1}^{t/\eta}}}{\sqrt{\frac{\vbb (t)}{1 - \beta_{2}^{t/\eta}}} + \epsilon} , \\
        \frac{\dd \mb}{\dd t} (t) &=  \frac{1}{\eta} (1 - \beta_1) \left ( \frac{\partial}{\partial \ab}  z (\ab (t) ) - \mb(t) \right ), \\
        \frac{\dd \vbb}{\dd t} (t) &= \frac{1}{\eta} (1-\beta_2) \left ( \left ( \frac{\partial}{\partial \ab} z (\ab (t) ) \right ) ^{2} - \vbb (t) \right ) .
        \label{eq:adam_ode_with_delta_t}
    \end{aligned}
\end{equation}
In Appendix \ref{appendix:adam_ode_verification}, we show how the discretization of Eq. \eqref{eq:adam_ode_with_delta_t} returns the original ADAM update rules in Eq. \eqref{eq:discrete_adam:param_update}-\eqref{eq:discrete_adam:unbiased_update}.

While we could attempt to learn each of the right-hand sides in Eq. \eqref{eq:adam_ode_with_delta_t}, this is costly and even unnecessary.
This would require building a library with $\ab$, $\mb$, and $\vbb$ (tripling the size of the state variables) and identifying an explicit dependence on $t$.
In the partially known case, compared to the expensive evaluation of $\partial z / \partial \ab$, the remainder of Eq. \eqref{eq:adam_ode_with_delta_t} is quite cheap.
Thus, we assume a partially known ODE, where we only seek to model the missing $\partial z / \partial \ab$ term with SINDy (which is also used within the variance update as $(\partial z / \partial \ab) ^ 2$).
Some recent extensions of SINDy focus on such discrepancy modeling \cite{ebers_discrepancy_2023, kaheman_learning_2019}.

Standard discrepancy modeling would only have access to observations generated by Eq. \eqref{eq:adam_ode_with_delta_t} and would need to isolate the unknown component, requiring the histories of $\mb$ and $\vbb$ (in addition to $\ab$).
However, in our case, we have direct access to observations of the unknown component $\pd{}{\mathbf{a}} z (\ab )$, simplifying this process considerably.
We simply record $\pd{}{\mathbf{a}} z (\ab )$ during the initial data-generating phase and build a SINDy model just for this expensive term.
While collecting $\Ab$ over iterations, we also record gradient observations:
$$
\Gb = \begin{bmatrix}
    \pd{}{\mathbf{a}} z (\ab (t_1) ) & \pd{}{\mathbf{a}} z (\ab (t_2) ) & \dots & \pd{}{\mathbf{a}} z (\ab (t_{\mathcal  K}) ) 
\end{bmatrix}^T \in \mathbb{R}^{\mathcal K \times n} .
$$
Then, we follow the same procedure as for the other methods, but with the modified goal
$$
\Gb \approx \mathbf{\Phi} ( \Ab ) \mathbf{\Xi} ,
$$
where we can reuse the same \texttt{STLSQ} solver as before to find $\mathbf{\Xi}$.
We replace $\partial z / \partial \mathbf{a}$ with the model $\fsindy (\ab) = \mathbf{\Xi}^T \boldsymbol{\phi} (\ab (t)) $, yielding a new ODE system:
\begin{equation}
    \begin{aligned}
        \frac{\dd \ab}{\dd t} (t) &= - \frac{\frac{\mb (t)}{1 - \beta_{1}^{t/\eta}}}{\sqrt{\frac{\vbb (t)}{1 - \beta_{2}^{t/\eta}}} + \epsilon} , \\
        \frac{\dd \mb}{\dd t} (t) &=  \frac{1}{\eta} (1 - \beta_1) ( \fsindy (\ab (t) ) - \mb(t) ), \\
        \frac{\dd \vbb}{\dd t} (t) &= \frac{1}{\eta} (1-\beta_2) ( \fsindy (\ab (t) ) ^{2} - \vbb (t) ) .
        \label{eq:adam_ode_with_delta_t_learned}
    \end{aligned}
\end{equation}
We integrate Eq. \eqref{eq:adam_ode_with_delta_t_learned} from $t =\mathcal  K \cdot \eta$ to $t = \mathcal M \cdot \eta$, beginning from $\ab (t_{\mathcal K})$, $\mb (t_{\mathcal K})$, and $\vbb (t_{\mathcal K})$.
The evolved optimization variables $\ab (t_{\mathcal M})$ should approximate the result of performing ADAM with the expensive objective gradient after $\mathcal M$ total iterations.
This method can be used with or without mini-batching of the objective function, and we showcase both in a subsequent numerical example.

While the interpretation of optimization as a dynamical system is an important motivation, it is not absolutely necessary for this method to function.
As we see in Eq. \eqref{eq:adam_ode_with_delta_t_learned}, the core requirement is the availability of an accurate surrogate for the optimizer’s update direction.
In principle, our method can be applied to other optimization methods that do not admit time continuous ODE forms.
We reiterate that, nonetheless, in some cases, the continuous perspective can provide some benefit \cite{flouris_explicit_2025, jacot_neural_2018, jentzen_convergence_2023} and also acts as a bridge to data-driven dynamics literature.



\subsection{Benefits and limitations of surrogate training dynamics}

After an initial data-generating phase, the SINDy models presented in Sec. \ref{sec:learned_gradient_flow}, \ref{sec:learned_newtons}, and \ref{sec:learned_adam} avoid any evaluation of the objective function $z$ or its derivatives.
This has both substantial benefits and possible consequences.
Most obviously, evaluating any terms related to $z$ is computationally expensive.
In engineering applications, this often corresponds to numerically simulating a system of interest, such as by solving partial differential equations (PDEs).
In machine learning, computing the loss and its gradient entails forward and backward passes through large models over substantial datasets, making it \textit{the} most costly part of training.
Thus, the primary aim of our methods is to avoid these costly computations by instead relying on the discovered training dynamics.

The utility of using a surrogate model for the evolution of the optimization variables fundamentally depends on its generalization. Because the learned dynamics are based on a finite set of observations, the performance of the surrogate hinges on how well the data-generating phase represents the training as a whole, and the evolution of the SINDy model provides no direct measure of the training quality (such as the loss). To address these concerns, we propose a heuristic variant of the general procedure outlined in Alg. \ref{alg:learned_gradient_flow}, focusing on the case of gradient descent, without loss of generality to Newton, ADAM, or other base optimizers.

\subsection{Scheduled retraining}
\label{sec:methods:retraining}

\begin{figure}[htb]
    \centering
    \includegraphics[width=\linewidth]{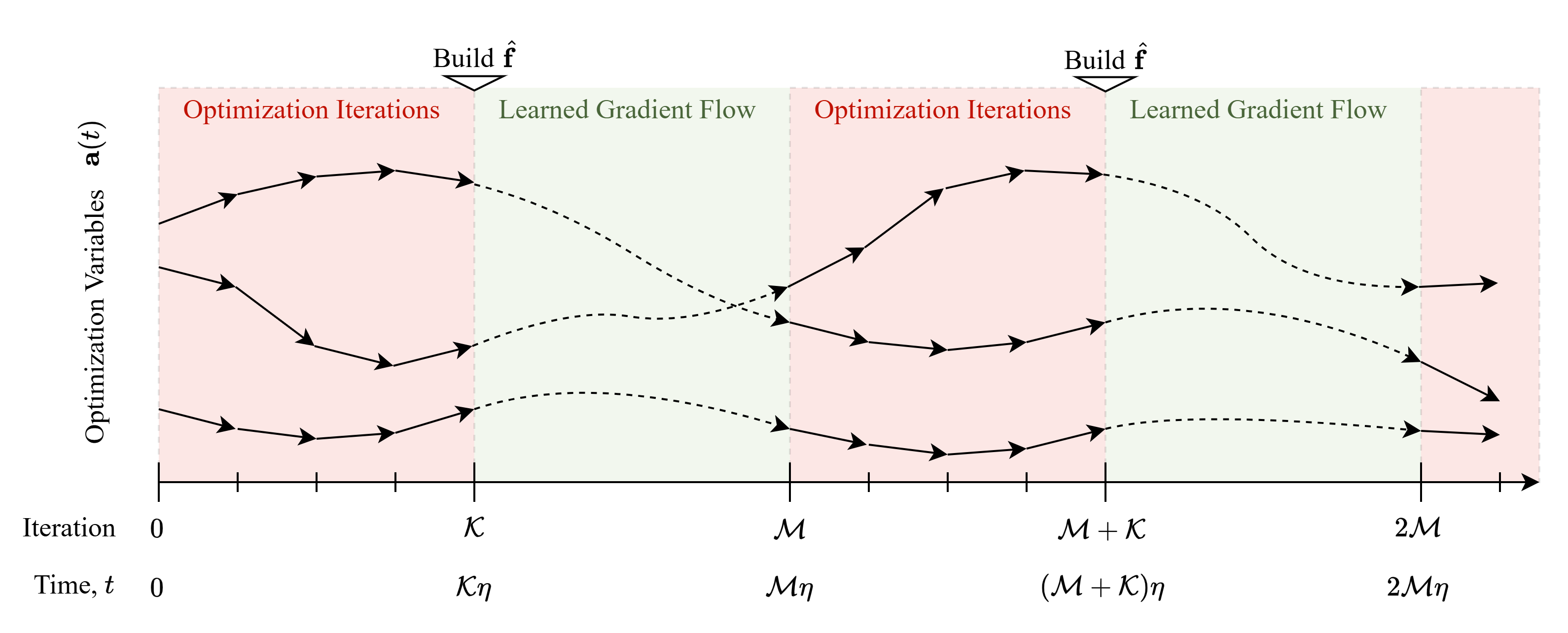}
    \caption{Scheduled retraining of gradient flow. We alternate between collecting histories for $\mathcal K$ iterations from the true optimizer with expensive gradient computations (shaded red) and applying the learned gradient flow, ending at iteration $\mathcal M$ (shaded green). Transitioning between these two steps, we must build $\fsindy$ as specified in Sec. \ref{sec:learned_gradient_flow}. The lines in the plot indicate trajectories of individual optimization variables.}
    \label{fig:scheduled_retraining_schematic}
\end{figure}

A preliminary approach assumes that these surrogate models decay in accuracy with time. In other words, a model is accurate at times close to the data-generating phase and less accurate when applied at later times. Under this assumption, we present a simple way to improve the generalizability of our method by devising an outer loop around Alg. \ref{alg:learned_gradient_flow} where we alternate between collecting optimization variable histories and applying a surrogate model. For each iteration of this outer loop, we iterate $\mathcal  K$ steps through the original optimizer to build the history, construct the surrogate model for the dynamics, integrate in time to reach iteration $\mathcal M$, and start over from the current values of the optimization variables.
This is equivalent to walking through Alg. \ref{alg:learned_gradient_flow} once, passing the output $\ab (t_{\mathcal M})$ as the initial optimization variables $\ab_1$ to another application of Alg. \ref{alg:learned_gradient_flow}, and so on. Adjusting the ratio between the history size $\mathcal K$ and the retrain interval $\mathcal M$ changes how much we rely on some model of the dynamics and thus the savings on the gradient evaluations.
For instance, setting $\mathcal M = \mathcal K$ results in always evaluating the true gradient and applying the base optimizer, while taking a large $\mathcal M$ means applying the dynamics model for a long duration.
Fig. \ref{fig:scheduled_retraining_schematic} intuitively illustrates this concept of alternating between training LGF and applying it.

\subsection{High-dimensional optimization problems}

For a state $\ab \in \mathbb{R}^n$ with $\mathcal K$ recorded iterations, the library matrix is $\mathbf{\Phi} \in \mathbb{R}^{\mathcal K \times p}$, where according to Eq. \eqref{eq:polynomial_size_scaling} $p$ scales based on $n$ depending on the polynomial order (e.g., $p = \mathcal{O}(n)$ for linear library terms and $\mathcal{O}(n^2)$ for quadratic library terms).
In many applications, the optimization variables are very high-dimensional, i.e., $n \gg 1$, making $p$ possibly much larger.
Storing and using this matrix $\mathbf{\Phi}$ becomes infeasible in such cases.
Thus, rather than building our dynamics model in $\mathbb{R}^n$, we will project to a smaller latent space $\mathbb{R}^r$, model the dynamics here, and then lift back to the original space.
For the measurement matrix $\Ab \in \mathbb{R}^{\mathcal K \times n}$, we can do this through the rank $r$ singular value decomposition (SVD):
$$
\mathbf{A}^T \approx \mathbf{U} \mathbf{\Sigma} \mathbf{V}^T ,
$$
assuming the state variables admit a low-rank representation, i.e., $r\ll n$. Here, $\mathbf{\Sigma} \in \mathbb{R}^{r \times r}$ is a diagonal matrix of $r$ largest (in descending order) singular values of $\mathbf{A}^T$, $\mathbf{U} \in \mathbb{R}^{n \times r}$ contains the corresponding $r$ left singular vectors (modes), and rows of $\mathbf{V} \in \mathbb{R}^{\mathcal K \times r}$ give the mode coefficients.
With this decomposition in hand, we can map $\ab$ to the latent space via $\abt = \mathbf{U}^T \ab$,
and recover $\ab$ with $\ab \approx \mathbf{U} \abt$. To create a model in the latent space, we must first project the optimization variables into the latent space via
$$
\Abt^T = \mathbf{U}^T \Ab^T \approx \mathbf{\Sigma} \mathbf{V}^T .
$$
We then continue with the process outlined in Sec. \ref{sec:sindy_for_ode_discovery}, but now on this projected space.
We define a candidate library function $\tilde{\boldsymbol{\phi}}: \mathbb{R}^r \to \mathbb{R}^{\tilde{p}}$ as all polynomials of up to some total order $\mathcal P$ (including cross terms), where $\tilde{p}$ is the number of these functions based on $r$.
Then, repeatedly applying this library function to the $\mathcal K$ rows of the projected optimization history gives the projected library matrix $\tilde{\mathbf{\Phi}} (\Abt) \in \mathbb{R}^{\mathcal K \times \tilde{p}}$.
Finally, the coefficients of the dynamics in the latent space $\tilde{\mathbf{\Xi}} \in \mathbb{R}^{\tilde{p} \times r}$ are computed such that
$$
\dot{\Abt} \approx \tilde{\mathbf{\Phi}} (\Abt) \tilde{\mathbf{\Xi}} .
$$
Here, $\dot{\Abt}$ can be computed as usual, but from $\Abt$ (e.g. via finite differencing), or alternatively via $\dot{\Abt}^T = \mathbf{U}^T \dot{\Ab}^T$.
In the case of recorded gradients $\Gb$ as in Sec. \ref{sec:learned_adam}, we compute
$$
\tilde{\Gb}^T = \mathbf{U}^T \Gb^T ,
$$
to build the latent model
$$
\tilde{\Gb} \approx \tilde{\mathbf{\Phi}} (\Abt) \tilde{\mathbf{\Xi}} .
$$
After finding $\tilde{\mathbf{\Xi}} \in \mathbb{R}^{\tilde{p} \times r}$, we can combine the dimensionality reduction with the dynamics model to prescribe a cheaper dynamics model on the original space:
$$
\frac{\dd \ab}{\dd t} (t) = \mathbf{U} 
    \, \tilde{\mathbf{\Xi}}^T \tilde{\boldsymbol{\phi}} (\mathbf{U}^T \ab (t)).
$$
This allows us to handle high-dimensional, $n \gg 1$, optimization problems.

Nonetheless, high-dimensional optimization problems may present difficulties due to storing and operating on the history matrix $\mathbf{A} \in \mathbb{R}^{\mathcal K \times n}$, yet we remark that this seemingly large memory requirement is present elsewhere in modern optimization.
Namely, the limited memory version of BFGS (LBFGS) \cite{nocedal_numerical_2006} has a comparable memory requirement.
For the LBFGS history size $\mathcal K'$, LBFGS stores roughly $2\mathcal K'$ vectors each of size $n$.
Anecdotally, in Sec. \ref{sec:results:top_op_dim_reduc}, we present a high-dimensional problem and use $\mathcal K=10$, as compared to the default history size for LBFGS $\mathcal K' = 100$ \cite{virtanen_scipy_2020}.
In this example, we use roughly $20$ times less memory than the default LBFGS.
Of course, this is not a fair comparison as LBFGS may reach an optimum with fewer epochs due to its quasi-Newton nature, but it goes to show that other widely-used algorithms have comparable memory footprints.
Further, there are several dimensionality reduction methods that have smaller memory requirements and accomplish the same purpose, including various randomized sketching approaches, which may require only a single pass \cite{dunton_pass-efficient_2020, dunton_deterministic_2021}.

\subsection{Summary of hyperparameters}

Before empirically demonstrating the methods presented in this section, we provide a brief overview of the important hyperparameters presented thus far. We remark that choosing the history size $\mathcal K$ and retrain interval $\mathcal M$ is a matter of both trial and error and experience with the optimization problem under study. Our numerical examples illustrate that different optimization problems permit different levels of acceleration with the surrogate model.

\begin{table}[h!]
    \centering
    \renewcommand{\arraystretch}{1.3}
    \begin{tabular}{c|p{10cm}}
        Hyperparameter & \multicolumn{1}{c}{Description} \\ \hline
        Learning rate, $\eta$ &
        The learning rate used by the base optimizer; also determines how frequently observations of the state are recorded or predicted \\
        History size, $\mathcal K$ & The number of samples of the state used to build $\fsindy$ \\
        Retrain interval, $\mathcal M$ & 
        The number of epochs specifying how often to discard the current $\fsindy$ and begin collecting data for building a new surrogate (which are counted towards the count, thus $\mathcal M \geq \mathcal K$); this determines the total savings \\
        Truncation rank, $r$ &
        If applicable, the dimension of the latent space in which we build the latent dynamics model after doing dimensionality reduction, $r \leq n$ \\
        Epochs &
        The total number of steps taken between both the optimizer and surrogate dynamics model(s) over the whole training procedure \\
        Polynomial order, $\mathcal P$ &
        The highest (total) degree of the polynomials used as library functions by SINDy \\
    \end{tabular}
\end{table}



\section{Numerical Examples}
\label{sec:results}


We now provide five numerical examples from engineering mechanics to showcase the LGF optimizer and compare its performance to traditional optimization methods. The first is a dynamic inverse problem with two unknown material parameters, which is chosen as a first case because it allows us to visualize movement through the loss landscape. The second problem is density-based topology optimization of a linear elastic structure, which is very high-dimensional and thus necessitates using dimensionality reduction techniques to learn the optimization dynamics. The third problem is a forward solve for nonlinear heat conduction with a spectral basis, and it is used to test the learning of continuous-time dynamics of Newton-based optimization. The fourth problem is that of ``full wave inversion,'' another dynamic inverse problem which is known to be challenging due to the non-convexity of the loss. As such, we use this problem to test the learned surrogate for ADAM optimization, which is more robust to non-convexity. Finally, in the fifth example, we experiment with LGF on a stochastic optimization problem from physics-informed machine learning. In particular, we solve for the displacement field of a three-dimensional linear elastic structure discretized by a neural network, where the objective function is computed with Monte Carlo integration. This example explores the ability of LGF to perform in the setting of noisy gradients. In each example, we compare LGF to a standard optimizer used to solve the problem at hand. Hyperparameter settings for the LGF optimizer are specified on a case-by-case basis, as the following examples show.

\subsection{Inverse problem with gradient descent}
\label{sec:eresults:ex1_2var}

As an initial illustration of our method, we consider a two-parameter inverse problem to visualize how the LGF optimizer approximates the true gradient flow dynamics. In this problem, we recover the parameterized conductivity of a heated bar from measurements. The space-time evolution of the temperature $u(x,t)$ in a bar of unit length is governed by 
\begin{equation}\label{ex1_eq}
    \pd{u}{t} = \pd{}{x}\qty( \kappa(x;a_1,a_2) \pd{u}{x} ) + b(x,t),\quad u(0,t) = u(1,t) = 0, \quad u(x,0) = 0,
\end{equation}
where $b(x,t)=2000\sin(2\pi x) \sin(2\pi t)$ is a prescribed heat source and the parameterized conductivity $\kappa(x;a_1,a_2)$ defines a multi-material bar as 
\begin{equation*}
    \kappa(x;a_1,a_2) = \begin{cases} a_1, \quad 0 \leq x \leq 0.5,\\
    a_2, \quad 0.5 < x \leq 1.
    \end{cases}
\end{equation*}

Note that the variable $t$ here is physical time rather than the pseudo-time of the optimization dynamics as in Eq. \eqref{eq:gradient_flow}. With that in mind, the initial temperature distribution is taken to be zero, and the bar has homogeneous Dirichlet boundaries on both ends. The goal of the inverse problem is to estimate the material parameters $a_1$ and $a_2$ from measurement data. To this end, we construct a numerical solution of Eq. \eqref{ex1_eq}.  The temperature is discretized with
\begin{equation}\label{ex1_dis}
    u(x,t) \approx \sum_{i=1}^N u_i(t) \sin( i \pi x),
\end{equation}
which ensures that the boundary conditions are satisfied automatically. Plugging in Eq. \eqref{ex1_dis}, the Bubnov-Galerkin weak form of Eq. \eqref{ex1_eq} is given by 
\begin{equation}\label{ex1_weak}
    \sum_{i=1}^N \pd{u_i}{t} \int_0^1 \sin( i \pi x) \sin( j \pi x) dx = -\sum_{i=1}^N u_i \int_0^1 \kappa(x;a_1,a_2) ij \pi^2 \cos(i \pi x) \cos(j \pi x) dx + \int_0^1 b(x,t) \sin( j \pi x) dx,
\end{equation}
where we have integrated by parts to transfer a spatial derivative onto the test function in the second term on the right-hand side. Defining the material parameters $\mathbf{a} = [a_1,a_2]^T$, Eq. \eqref{ex1_weak} can be written more compactly as
\begin{equation}\label{ex1_discrete}
    \mathbf{M} \pd{\mathbf u}{t} + \mathbf{K}(\mathbf{a}) \mathbf u = \mathbf{F}(t), \quad \mathbf u(0) = \mathbf{0}.
\end{equation}

The mass matrix $\mathbf{M}$, the stiffness matrix $\mathbf{K}$, and the force vector $\mathbf{F}(t)$ are defined through the indexed integral quantities of Eq. \eqref{ex1_weak}. Assume that we have full-field measurement data $v(x,t)$ over the interval $t\in[0,T]$, which is taken from the system without noise for a given $\mathbf{a}$. The inverse problem is then
\begin{equation}\label{ex1_obj}
\begin{aligned}
    &\underset{\boldsymbol a, \mathbf u}{\text{argmin }} \frac{1}{2} \int_0^T\int_0^1\Big(  \sum_{i=1}^N u_i(t) \sin( i \pi x)- v(x,t)\Big)^2 dx dt \\
    & \text{s.t. }\mathbf{M} \pd{\mathbf u}{t} + \mathbf{K}(\mathbf{a}) \mathbf u - \mathbf{F}(t) = 0, \quad \mathbf u(0) = \mathbf 0.
\end{aligned}
\end{equation}

If we treat the temperature parameters $\mathbf u$ as an explicit function of the material parameters $\boldsymbol a$ through the governing equation, the inverse problem can be written without constraints. Using the backward Euler method to time integrate, the temperature parameters are obtained iteratively as a function of the material parameters with
%
\begin{equation}\label{ex1_time}
   \qty( \frac{\mathbf{M}}{\Delta t} + \mathbf{K}(\mathbf{a}) ) \mathbf u_{k+1} = \qty( \mathbf{F}_{k+1} + \frac{\mathbf{M}}{\Delta t} \mathbf {u}_k).
\end{equation}
We write temperature parameters which are computed through Eq. \eqref{ex1_time} as $\mathbf u(t;\mathbf{a})$. The unconstrained form of the inverse problem is then
\begin{equation}\label{ex1_obj_reduced}
\begin{aligned}
    & \underset{\mathbf{a} }{\text{argmin }} z(\mathbf{a}),\\
    & z(\mathbf{a}) = \frac{1}{2} \int_0^T\int_0^1\Big(  \sum_{i=1}^N u_i(t; \mathbf{a}) \sin( i \pi x)- v(x,t)\Big)^2 \dd x \, \dd t.
\end{aligned}
\end{equation}

\begin{figure}[hbt!]
\centering
\includegraphics[width=0.99\textwidth]{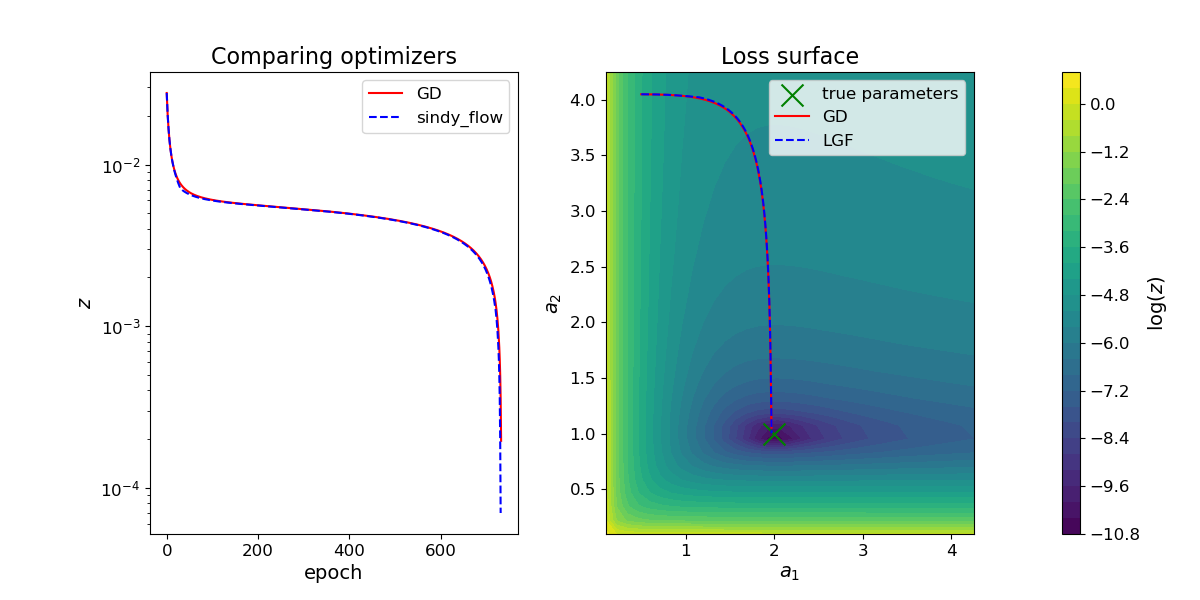}
\caption{Convergence of the two optimizers to the minimum of the data loss defining the inverse problem. The optimizer based on the SINDy surrogate for the dynamics tracks closely with standard gradient descent and accelerates the optimization by $200 \%$.}
\label{ex1_sf}
\end{figure}

The measurement data $v(x,t)$ is generated by time integrating Eq. \eqref{ex1_discrete} with the backward Euler method and $N=30$ basis functions in the spatial discretization. The true material parameters are $\mathbf{a}_{\text{true}} = [2,1]^T$. In solving Eq. \eqref{ex1_obj_reduced}, we compare standard gradient descent (GD) against the LGF optimizer with a linear approximation of the continuous-time optimization dynamics. The history size id $\mathcal K=10$, and the retrain interval is $\mathcal M=30$. The gradient $\partial z / \partial \mathbf{a}$ is computed using automatic differentiation with PyTorch by backpropagating through the time integration scheme of Eq. \eqref{ex1_time}. With a fixed number of epochs, we note that our chosen parameter settings for the surrogate model indicate $3\times$ savings in terms of gradient evaluations. The learning rate for both optimizers is set at $\texttt{lr}=1 \times 10^{-2}$, and the optimization is run for $700$ epochs, which we observe to be sufficient for convergence. See Table \ref{tab:ex1} for a summary of the problem setup and Figure \ref{ex1_sf} for the results of the comparison. We define the acceleration of the method as $100 \times (\mathcal M / \mathcal K - 1)$. The recovered parameters from the two methods are 
\begin{equation*}
    \mathbf{a}_{\text{GD}} = \begin{bmatrix}
        1.97 \\ 1.06
    \end{bmatrix}, \quad \mathbf{a}_{\text{LGF}} = \begin{bmatrix}
        1.97 \\ 1.04
    \end{bmatrix}.
\end{equation*}


Both achieve similar accuracies in recovering the parameters in the dynamic inverse problem and have similar final losses ($z(\ab_{\text{GD}}) = 2 \times 10^{-4}$ and $z(\ab_{\text{LGF}}) = 1 \times 10^{-4}$), though the LGF model saves $20$ out of every $30$ evaluations of the objective gradient, which involves back-propagating through the time integration of the temperature field.

\begin{table}
    \centering
    \begin{tabular}{|c|c|}
        \hline
        \textbf{Problem parameter} &  \textbf{Setting} \\
        \hline
        $\#$ optimization variables & $2$  \\
        $\texttt{lr} $& $1 \times 10^{-2}$  \\
        Retrain interval, $\mathcal M$ & $30$ \\
        History size, $\mathcal K$ & $10$\\
        Polynomial order, $\mathcal P$ & $1$ \\
        $\texttt{epochs}$ & $700$ \\
        Truncation rank, $r$ & N/A \\
         \hline
    \end{tabular}
    \caption{Problem parameters for the 1D dynamic inverse problem.}
    \label{tab:ex1}
\end{table}

\subsection{Topology optimization with gradient descent and dimensionality reduction}
\label{sec:results:top_op_dim_reduc}

We now test our method on a standard high-dimensional optimization problem from engineering mechanics: density-based topology optimization with a strain energy objective. In continuous form, this optimization problem reads
\begin{equation}\label{to_continuous}
\begin{aligned}
    & \underset{ \rho(\mathbf{x}), \mathbf{u}(\mathbf{x})}{\text{argmin }} \int_{\Omega} \frac{1}{2} \epsilon_{ij} \rho(\mathbf{x})^pC_{ijk\ell} \epsilon_{k\ell} \, \dd \Omega \\
    & \text{s.t. } \pd{\sigma_{ij}}{x_j} + b_i = 0, \quad  \int \rho(\mathbf{x}) \, \dd \Omega - V_0 = 0, \quad  0 \leq \rho \leq 1,
\end{aligned}
\end{equation}
where $C_{ijk\ell}$ is the linear elastic constitutive tensor for an isotropic, homogeneous material, $p$ is a ``SIMP'' penalty parameter \cite{bendsoe_topology_2004}, $V_0$ is the volume constraint, $\rho(\mathbf{x})$ is the density field, and the strain and stress are defined in terms of the displacement $\mathbf{u}(\mathbf{x})$ as 
\begin{equation*}
    \epsilon_{ij} = \frac{1}{2}\qty( \pd{u_i}{x_j} + \pd{u_j}{x_i}), \quad \sigma_{ij} = \rho^pC_{ijk\ell} \epsilon_{k\ell}.
\end{equation*}

Working with a finite element discretization, the displacement components are stored on the nodes of the finite element mesh, and the density is discretized as constant over elements. The governing equation for static equilibrium in Eq. \eqref{to_continuous} is converted to the Bubnov-Galerkin weak form. Furthermore, we can enforce the box constraint on the density automatically by defining the elemental densities through the optimization variables $\mathbf{a}$. The densities are then computed as
\begin{equation}\label{density}
    \boldsymbol \rho = \frac{1}{2} \Big( \tanh( \mathbf{a}) + 1\Big),
\end{equation}
The discretized optimization problem can be written as
\begin{equation}\label{discrete_const}
\begin{aligned}
    & \underset{\mathbf{a}, \mathbf{u}}{\text{argmin }} \frac{1}{2} \mathbf{u}^T \mathbf{K}(\mathbf{a}) \mathbf{u} \\
    & \text{s.t. } \mathbf{K}(\mathbf{a}) \mathbf{u} - \mathbf{F} = \mathbf{0}, \quad \sum_{e=1}^{n_{el}} \rho( a_e) \Omega_e - V_0 = 0,
\end{aligned}
\end{equation}
where $n_{el}$ is the number of elements in the finite element mesh, $\Omega_e$ are the volumes of each element, $\mathbf{u}$ are the nodal displacements, $\mathbf{K}$ is the stiffness matrix, and $\mathbf{F}$ is the force vector. We convert this problem to an unconstrained form by treating the displacement as an explicit function of the optimization variables through the governing equation for equilibrium and by transforming the volume constraint into a penalty. At each setting of the optimization variables $\mathbf{a}$, the displacement is computed through the governing equation of equilibrium with
%
\begin{equation}
\label{eqn:equi}
    \mathbf{K}(\mathbf{a}) \mathbf{u}(\mathbf{a}) = \mathbf{F}.
\end{equation}

As such, the unconstrained form of the topology optimization problem of Eq. \eqref{discrete_const} is 
\begin{equation}\label{discrete_unconst}
\begin{aligned}
    & \underset{\mathbf{a}}{\text{argmin }} z(\mathbf{a}) ,\\
    & z(\mathbf{a}) = \frac{1}{2} \mathbf{u}(\mathbf{a})^T \mathbf{K}(\mathbf{a}) \mathbf{u}(\mathbf{a}) + \frac{\beta}{2}\Big( \sum_{e=1}^{n_{el}} \rho(a_e) \Omega_e- V_0 \Big)^2,
\end{aligned}
\end{equation}
where $\beta$ is a penalty hyperparameter. We can enforce that the displacement satisfies the governing equation at each design iteration in the computation of the gradient of this loss. Differentiating the objective in Eq. \eqref{discrete_unconst}, we note that both the stiffness matrix and the displacement depend on the optimization variables. The gradient is thus
\begin{equation}\label{grad_obj}
  \pd{z}{a_j} = \frac{1}{2} \mathbf{u}^T \pd{\mathbf{K}}{a_j} \mathbf{u} + \mathbf{u}^T\mathbf{K} \pd{\mathbf{u}}{a_j} + \beta\Big( \sum_{e=1}^{n_{el}} \rho(a_e) \Omega_e - V_0 \Big) \frac{\Omega_j}{2}(1-\tanh^2(a_j) ).
\end{equation}

We can obtain the ``sensitivity derivative'' $\partial \mathbf{u} / \partial \mathbf{a}$ by differentiating the governing equation Eq. \eqref{eqn:equi}:
\begin{equation}\label{sens}
    \mathbf{K}(\mathbf{a})\pd{\mathbf{u}}{\mathbf{a}} = -\pd{\mathbf{K}}{\mathbf{a}}\mathbf{u}(\mathbf{a}).
\end{equation}
Plugging Eq. \eqref{sens} into Eq. \eqref{grad_obj}, the gradient of the objective takes on a particularly simple form
\begin{equation}\label{to_grad}
        \pd{z}{\mathbf{a}} = -\frac{1}{2}\mathbf{u}^T \pd{\mathbf{K}}{\mathbf{a}} \mathbf{u} + \beta \Big( \sum_{e=1}^{n_{el}} \rho(a_e) \Omega_e - V_0 \Big) \frac{\boldsymbol \Omega}{2}(1-\tanh^2(\mathbf{a}) ),
\end{equation}
where $\boldsymbol \Omega$ is the vector of element volumes. The strain energy gradient in Eq. \eqref{to_grad} can be assembled from elemental contributions. Using Eq. \eqref{density} and the SIMP penalization of the elemental stiffness matrix, the elemental contribution to the strain energy gradient is given as 
\begin{equation}\label{elemental}
    \qty(-\frac{1}{2}\mathbf{u}^T \pd{\mathbf{K}}{\mathbf{a}} \mathbf{u})_j = -\frac{1}{4}  p \rho(a_j)^{p-1}( 1- \tanh^2(a_j)) (\mathbf{u}_j^e)^T\mathbf{K}_e (\mathbf{u}_j^e) ,
\end{equation}
where $\mathbf{K}_e$ is the elemental stiffness matrix and $\mathbf{u}^e_j$ is the elemental displacement vector corresponding to the $j$-th element. With Eq. \eqref{elemental}, the gradient of the strain energy can be assembled by looping over elements in the finite element mesh.

\begin{figure}[hbt!]
\centering
\includegraphics[trim = 20mm 20mm 20mm 20mm, clip, width=0.65\textwidth]{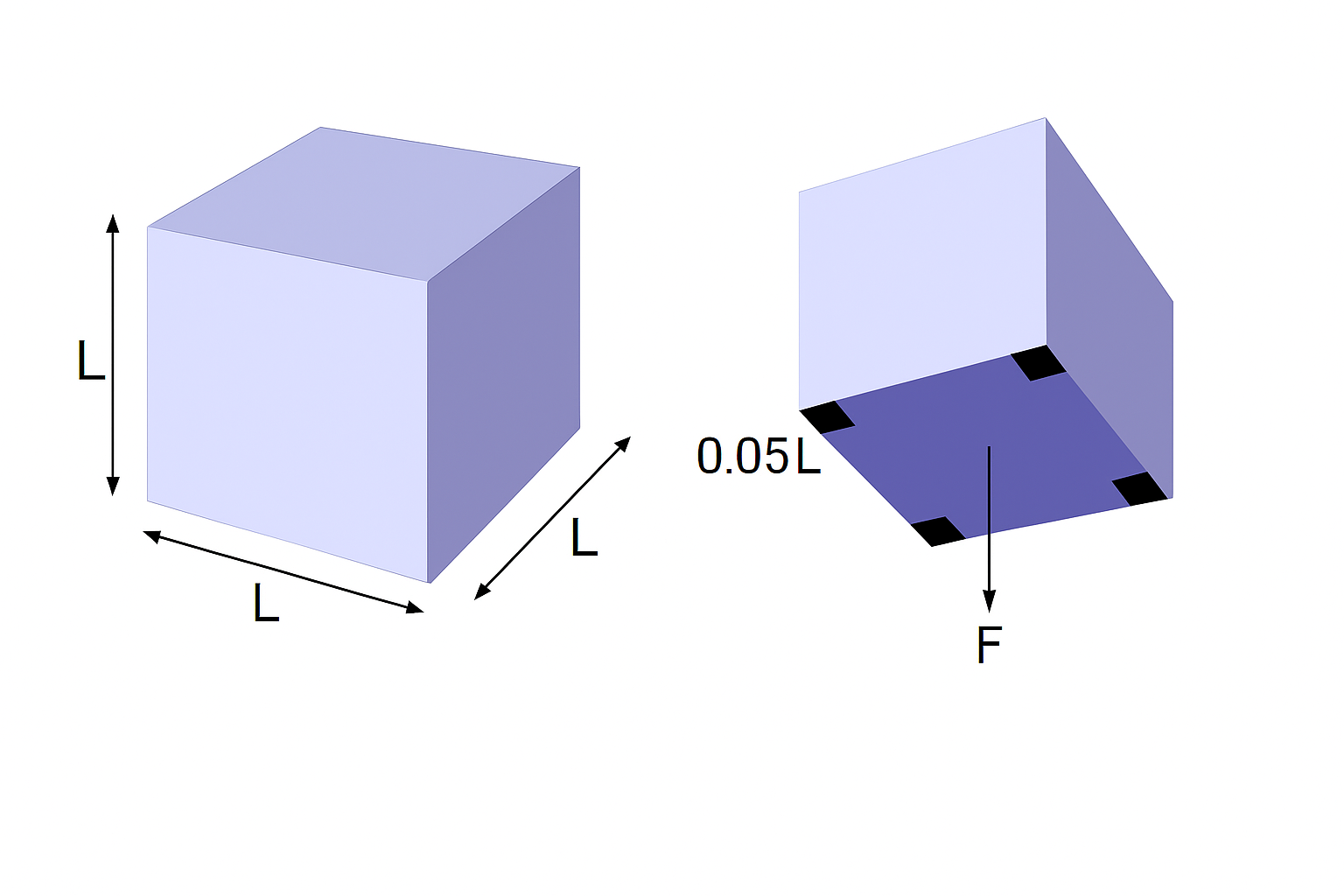}
\caption{Design domain, boundary conditions, and loading for the compliance minimization design problem. We take $L=1$ and use the symmetry of the problem to only simulate a quarter of the domain.}
\label{to_diagram}
\end{figure}

\begin{figure}[hbt!]
\centering
\includegraphics[trim = 5mm 10mm 20mm 15mm, clip, width=1\textwidth]{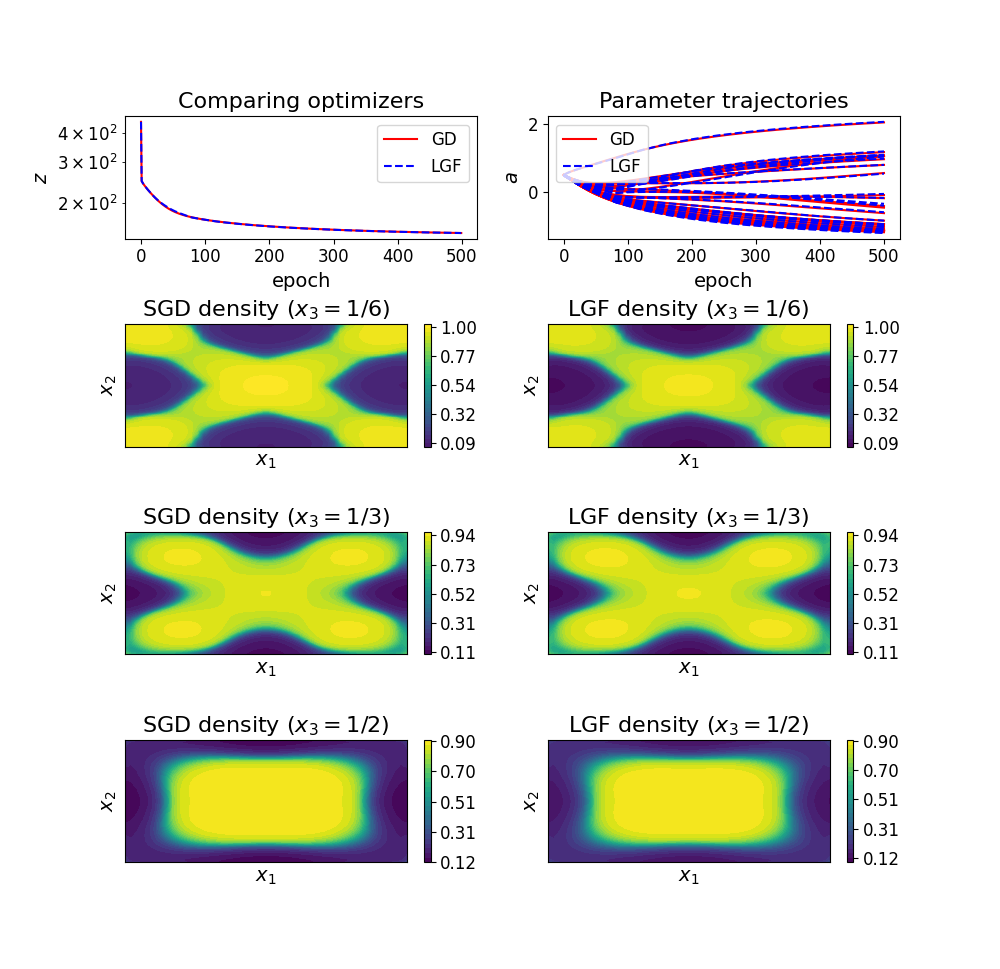}
\caption{Results of the density-based topology optimization with $100\%$ acceleration using the surrogate model. We build the linear surrogate in a $2$-dimensional subspace. This dramatic reduction in the dimensionality (from $n=13824$ to $r=2$) keeps surrogate evaluations lightweight once constructed, as each update is carried out in a rank-$r$ subspace.
}
\label{to_ex}
\end{figure}

With an unconstrained objective and its gradient defined, we can set up a problem and carry out the topology optimization. See Figure \ref{to_diagram} for the problem setup. We take the design domain to be the unit cube ($L=1$), but use symmetry to define the computational domain as $\Omega=[0,1/2] \times [0,1/2] \times [0,1]$. The displacement field is clamped in a cube with side length $0.05$ around the origin. We use 8-node brick finite elements with linear interpolations of the displacement components and $24$ elements per side of the computational domain. This amounts to $13824$ equal-sized elements and $46875$ degrees of freedom before applying the Dirichlet and symmetry boundary conditions. We compare the LGF optimizer to standard gradient descent. We take the learning rate to be $\texttt{lr}=1$, and aim for $100\%$ cost savings with history size $\mathcal{K} = 20$ and retrain interval $\mathcal{M} =40$. We use a linear surrogate for dynamics with polynomial order $\mathcal P=1$ and build the dynamics in a low-dimensional subspace with truncation rank $r=2$.
Building the reduced-order LGF surrogate requires computing a rank-$r$ basis, which has cost $\mathcal{O}(n\,r\,\mathcal{K})$ for truncated SVD.
This procedure is performed once per retraining interval $\mathcal{M}$, so its cost is
amortized over the subsequent surrogate steps.
By contrast, each ``true'' objective/gradient evaluation in topology optimization requires solving the sparse
equilibrium system $\mathbf{K}(\mathbf{a})\mathbf{u}=\mathbf{F}$ and assembling
sensitivities over elements.
Thus, while the computational cost of computing the truncated SVD (once per retrain interval) is not insignificant, it also does not dominate the cost of the problem by comparison to the regular training objective.
This amortized basis-construction cost mirrors the standard offline/online trade-off in reduced-order modeling,
where an upfront low-rank basis is built at nontrivial cost to enable many subsequent inexpensive evaluations.
Given the small dimension of the surrogate model, we do not enforce sparsity of the learned coefficients, setting both the ridge regression weight $\alpha$ and $\texttt{threshold}$ to $0$. See Figure \ref{to_ex} for the results. After $500$ iterations, the difference between the density fields is 
\begin{equation*}
    \frac{ \int_{\Omega} |\rho_{\text{GD}}(\mathbf{x}) - \rho_{\text{LGF}}(\mathbf{x}) |^2 d\Omega}{ \int_{\Omega}  | \rho_{\text{GD}}(\mathbf{x}) |^2 d\Omega} = 2.0 \times 10^{-4}.
\end{equation*}

The accuracy of the surrogate is striking not just because of the savings in objective/gradient evaluations, but especially because the dynamics are built in a two-dimensional subspace. This is an incredible contrast with the apparent dimensionality of the optimization problem of $13824$. Here, the dynamics of density-based topology optimization has an intrinsically low-dimensional and predictable dynamical structure, in spite of the complexity of the physics of the problem. We note that our choice to avoid box constraints by using hyperbolic tangent functions dramatically delays convergence of the density field to binary values, as Eq. \eqref{density} only converges to $0$ or $1$ in a limiting sense.

\begin{table}
    \centering
    \begin{tabular}{|c|c|}
        \hline
        \textbf{Problem parameter} &  \textbf{Setting} \\
        \hline
        $\#$ optimization variables & $13824$  \\
        Volume constraint ($V_0$) & $0.3$ \\
        SIMP penalty ($p$) & 4 \\
        $\texttt{lr} $& $1$  \\
        Retrain interval, $\mathcal M$ & $40$ \\
        History size, $\mathcal K$ & $20$\\
        Polynomial order, $\mathcal P$ & $1$ \\
        $\texttt{epochs}$ & $500$ \\
        Truncation rank, $r$ & $2$ \\
         \hline
    \end{tabular}
    \caption{Problem parameters for the 3D density-based topology optimization problem.}
    \label{tab:ex2}
\end{table}

\begin{figure}[hbt!]
\centering
\includegraphics[width=1\textwidth]{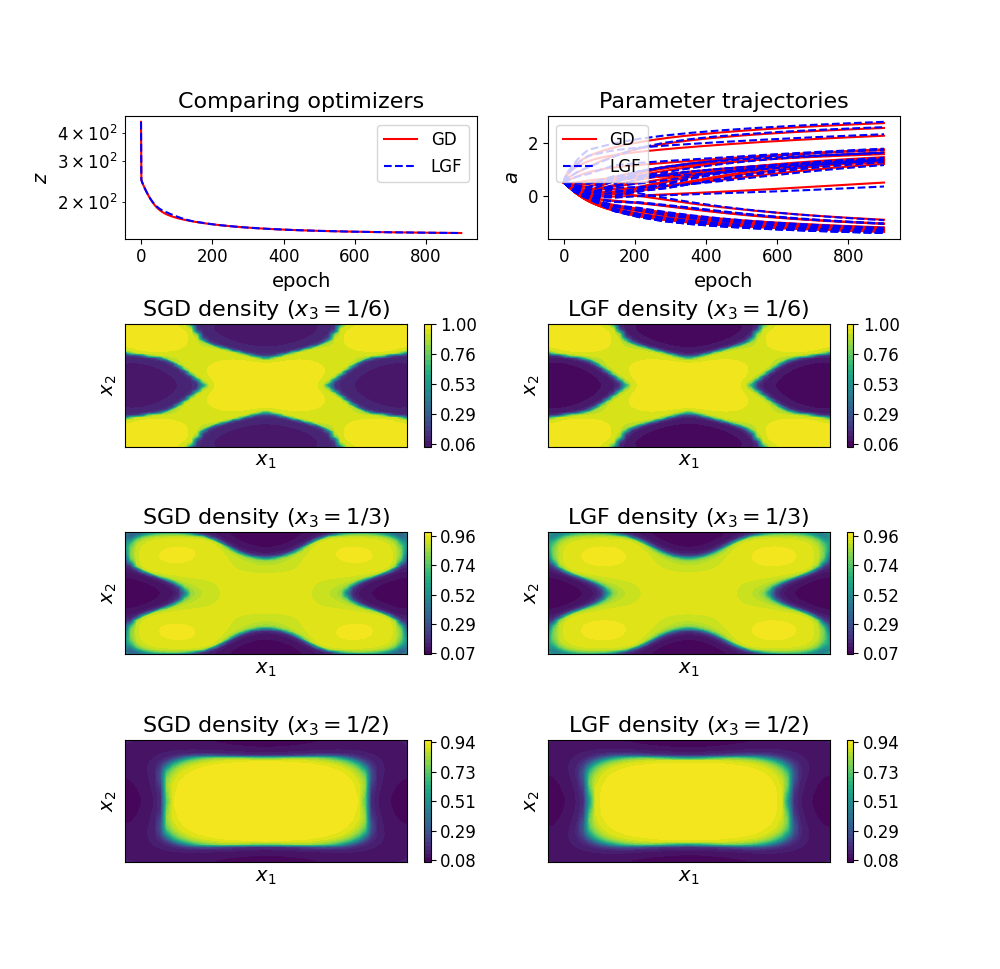}
\caption{Results of the density-based topology optimization with $250\%$ acceleration using the surrogate model. The learned model for the linear optimization dynamics is built in a $2$-dimensional subspace. The surrogate model again closely tracks the reference solution obtained with gradient descent.
}
\label{to_ex2}
\end{figure}

Keeping all other problem parameters constant, we run the problem again with a history size $\mathcal{K} = 20$ and retrain interval $\mathcal{M} =70$, corresponding to $250\%$ savings in evaluations of the objective and its gradient. The optimization is run for $900$ epochs. See Figure \ref{to_ex2} for the results. Even with the increased acceleration, the optimized design from the surrogate agrees with the design obtained with gradient descent. These examples suggest that, despite the complexity of the physics of the topology optimization problem, the optimization dynamics are simple. It is an interesting direction for future work to focus specifically on learned surrogates for topology optimization, and systematically study to what extent the optimization process can be accelerated before failure occurs.

\subsection{Nonlinear heat conduction with Newton's method}
\label{sec:ex:heat-newton}


In continuous time, Newton optimization is given by Eq. \eqref{newtonflow}. The LGF optimizer must learn the dynamics from historical data on the optimization trajectories, so, for LGF to be effective, an appropriate test problem must require a large number of Newton steps to converge. One such problem is the scalar $P$-Laplace equation---which is used as a nonlinear generalization of heat conduction \cite{lindqvist_nonlinear_nodate}---with the addition of radiative heat loss. Given its variational formulation, this problem allows us to conveniently represent a PDE forward solve as an optimization problem. In the continuous setting, the variational energy is given by

\begin{equation*}
    z( u(\mathbf{x})) = \int_{\Omega} \kappa(\mathbf{x}) | \nabla u \cdot \nabla u|^P + \frac{\sigma}{5}u^5 - b(\mathbf{x}) u \, \dd \Omega,
\end{equation*}
where $u(\mathbf{x})$ is the temperature field, $\kappa(\mathbf{x})$ is the conductivity, $P$ determines the order of the nonlinearity, $\sigma$ is the emissivity of the material, and $b(\mathbf{x})$ is a heat source. We note that choosing to raise $\nabla u \cdot \nabla u$ to the $P$-th power (as opposed to computing the $P$-th power of $|\nabla u|$) is non-standard but streamlines calculations and is inconsequential for this example. For our problem, we take the domain to be $\Omega = [0,1]^2$ and the boundaries to be homogeneous Dirichlet $u(\mathbf{x}) \Big |_{\partial \Omega} = 0$. We discretize the solution with a spectral basis per
\begin{equation*}
    u(\mathbf{x}) \approx \sum_{i=1}^M \sum_{j=1}^M  \tilde a_{ij} \sin( i \pi x_1) \sin( j \pi x_2)  := \sum_{i=1}^{M^2} a_i f_i(\mathbf{x}),
\end{equation*}
where the relation between the coefficients $\tilde a_{ij}$ and $a_i$ (as well as their corresponding basis functions) is a reshape operation. Plugging the discretization into the variational energy, the governing equation for the temperature field is
\begin{equation}\label{newton_opt}
\begin{aligned}
    & \underset{\mathbf{a}}{\text{argmin }} z(\mathbf{a})\\
    & z(\mathbf{a}) = \int_{\Omega} \qty[ \kappa(\mathbf{x}) \Big |  \sum_{i=1}^{M^2} \sum_{j=1}^{M^2} a_i a_j \nabla f_i \cdot \nabla f_j\Big|^P + \frac{\sigma}{5}\qty( \sum_{i=1}^{M^2} a_i f_i (\mathbf{x})  )^5 - b(\mathbf{x}) \sum_{i=1}^{M^2} a_i f_i(\mathbf{x}) ] d\Omega.
\end{aligned}
\end{equation}

\begin{figure}[hbt!]
\centering
\includegraphics[trim = 20mm 25mm 20mm 15mm, clip, width=1\textwidth]{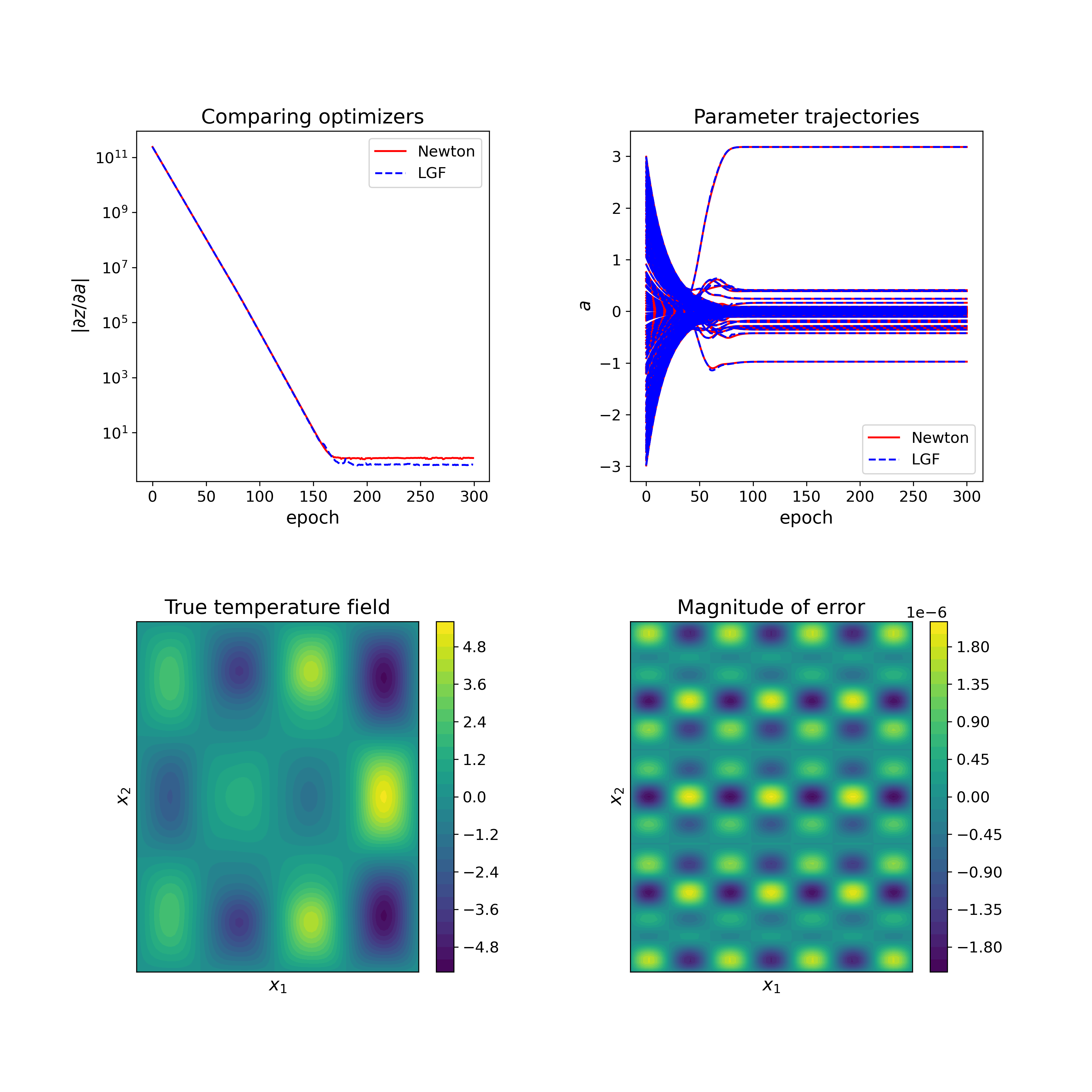}
\caption{Comparing the surrogate for the full Newton dynamics to the true dynamics. Using the LGF optimizer represents a $25\%$ acceleration. The surrogate briefly departs from the true Newton dynamics as the parameters begin to converge, but quickly recovers, leading to a final solution with comparable accuracy to the true Newton method.}
\label{newton_ex}
\end{figure}

Newton-type optimization requires a quadratic approximation of the objective, meaning that we need both the gradient of the objective and the Hessian matrix. These can be computed as 
\begin{equation*}
\begin{aligned}
    \pd{z}{a_k} = \int_{\Omega} \Bigg[  2P \kappa(\mathbf{x})\Big |  \sum_{i=1}^{M^2} \sum_{j=1}^{M^2} a_i a_j \nabla f_i \cdot \nabla f_j\Big|^{P-1} \sum_{j=1}^{M^2} a_j \nabla f_j \cdot \nabla f_k + \sigma \qty(\sum_{i=1}^{M^2} a_i f_i(\mathbf{x}) )^4 f_k(\mathbf{x}) - b(\mathbf{x}) f_k(\mathbf{x}) \Bigg] d\Omega  \\ 
    \frac{\partial^2 z}{\partial a_k \partial a_{\ell}} = \int_{\Omega} \Bigg[  4 \kappa(\mathbf{x}) P (P-1)\Big |  \sum_{i=1}^{M^2} \sum_{j=1}^{M^2} a_i a_j \nabla f_i \cdot \nabla f_j\Big|^{P-2} \qty(\sum_{j=1}^{M^2} a_j \nabla f_j \cdot \nabla f_k) \qty(\sum_{i=1}^{M^2} a_i \nabla f_i \cdot \nabla f_{\ell} ) \\ + 2P \kappa(\mathbf{x})\Big |  \sum_{i=1}^{M^2} \sum_{j=1}^{M^2} a_i a_j \nabla f_i \cdot \nabla f_j\Big|^{P-1} \nabla f_k \cdot \nabla f_{\ell} +  4\sigma \qty( \sum_{i=1}^{M^2} a_i f_i(\mathbf{x}))^3 f_k (\mathbf{x}) f_{\ell}(\mathbf{x}) \Bigg] d\Omega.
\end{aligned}
\end{equation*}
We will solve the Newton-type optimization problem for the temperature field given by Eq. \eqref{newton_opt} for a multi-material object. The problem parameters for the governing equation and discretization are as follows:
\begin{equation*}
\begin{aligned}
    & b(\mathbf{x}) = 10^{7}x_1\sin(4\pi x_1) \sin(3 \pi x_2 )  ,\\
    & \kappa(\mathbf{x}) = \begin{cases} 20, \quad \text{if } \text{max}( |\mathbf{x} - [1/2,1/2]^T |) \leq 0.25, \\
    1, \quad \text{else},
    \end{cases} \\
    & M = 15, \\
    & \mathbf{a}_0 \sim \mathcal{U}(-3,3),
\end{aligned}
\end{equation*}
where $\mathbf{a}_0$ is the initial guess for the coefficients. In the LGF optimizer, we relax the Newton solve with $\Delta t=\texttt{lr} = 0.15$ in order to ensure stability, take the polynomial order of the surrogate as $\mathcal P=1$, and use the history size $\mathcal K=15$ and retrain interval $\mathcal M=20$. Note that we use the same step relaxation of $0.15$ for the true Newton solve as well. We form the objective and its gradients with $5625$ equally spaced integration points throughout the unit square domain, and the optimization is run for 300 epochs. See Table \ref{tab:ex3} for a summary of the problem set up and Figure \ref{newton_ex} for the comparison of a standard Newton updating vs. the surrogate for the Newton dynamics. For both optimizers, the initial norm of the gradient is $|\partial z (\mathbf{a}_0)/ \partial \mathbf{a}| := R_0 = 2.38 \times 10^{11}$, and after 300 steps, the norm of the gradient is reduced to
\begin{equation*}
    R_{\text{NEWTON}} = 0.7, \quad R_{\text{LGF}} = 1.2.
\end{equation*}
Both residuals have been reduced by more than ten orders of magnitude. The normalized difference between the two temperature fields is 
\begin{equation*}
    \frac{ \int_{\Omega} | u_{\text{NEWTON}}(\mathbf{x}) - u_{\text{LGF}}(\mathbf{x}) |^2 d\Omega}{ \int_{\Omega}  | u_{\text{NEWTON}}(\mathbf{x}) |^2 d\Omega} = 1.7 \times 10^{-13},
\end{equation*}
which represents remarkable agreement between the Newton optimization and the surrogate for its dynamics, which saves $5$ out of every $20$ evaluations of the right-hand side of the continuous-time Newton dynamics. In the case of Newton optimization, evolving the dynamical system for optimization involves computing the gradient of the objective, its Hessian, the inverse of the Hessian, and the matrix-vector product of the inverse Hessian with the gradient. By comparison, the surrogate is extremely cheap to evaluate and leads to a solution that, in this case, is indistinguishably close to the true Newton method.

\begin{table}
    \centering
    \begin{tabular}{|c|c|}
        \hline
        \textbf{Problem parameter} &  \textbf{Setting} \\
        \hline
        $\#$ optimization variables & $225$  \\
        Emissivity ($\sigma$) & $4$ \\
        Order of nonlinearity ($P$) & $2$ \\
        $\texttt{lr} $& $1.5 \times 10^{-1}$  \\
        Retrain interval, $\mathcal M$ & $20$ \\
        History size, $\mathcal K$ & $15$\\
        Polynomial order, $\mathcal P$ & $1$ \\
        $\texttt{epochs}$ & $300$ \\
        Truncation rank, $r$ & N/A\\
         \hline
    \end{tabular}
    \caption{Problem parameters for the nonlinear heat conduction energy minimization problem.}
    \label{tab:ex3}
\end{table}

\begin{figure}[hbt!]
\centering
\includegraphics[trim = 20mm 25mm 20mm 15mm, clip, width=1\textwidth]{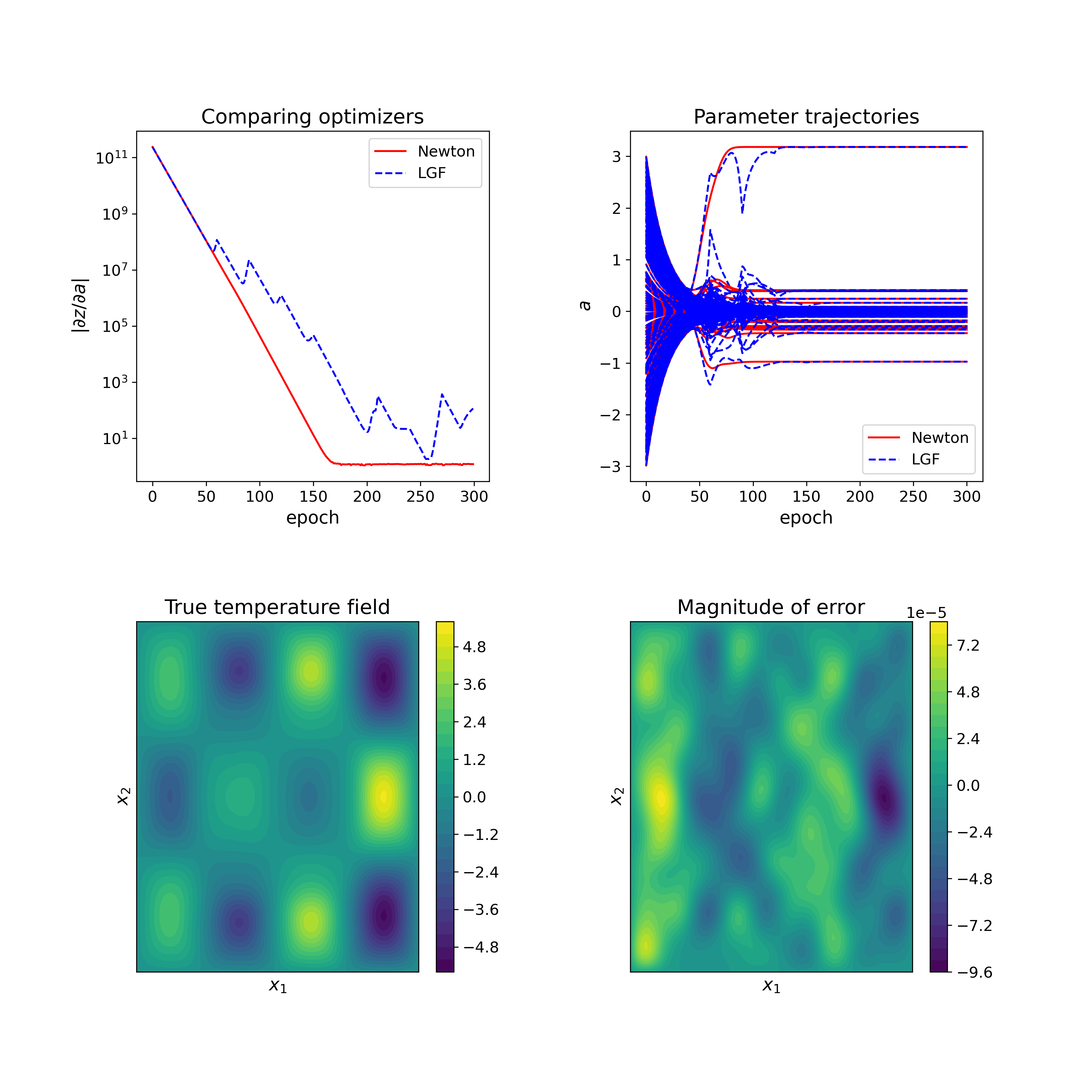}
\caption{Comparing the surrogate model against the true Newton dynamics with a $100\%$ acceleration. Because the acceleration is too aggressive, the surrogate model for the Newton dynamics ceases to be accurate and drives the norm of the residual up when active.}
\label{newton_ex2}
\end{figure}

Thus far, no guidelines have been provided on how to choose the acceleration of the learned gradient flow optimizer. As this is an initial work showcasing the possibility of accelerating optimization with equation discovery, future work is required to explore automating the acceleration. However, to provide some insight into the pitfalls of choosing an overly aggressive acceleration ratio, we re-run the problem with $\mathcal M=30$ while keeping all other problem parameters constant. See Figure \ref{newton_ex2} for the results. The two optimizers closely track each other in the initial stages of the optimization, yet the surrogate for the Newton dynamics consistently drives the norm of the residual up after epoch $175$. Though the magnitude of the error between the two solutions is still small, the surrogate model is taking steps that later need to be corrected by the true optimization dynamics. Here, because the aggressive acceleration consistently increases the loss, the surrogate provides little benefit. We note from experience that, as expected, even larger acceleration values correspond to increased discrepancy between the true and surrogate Newton dynamics.

\subsection{Full wave inversion with ADAM optimization}
\label{sec:ex:wave-inversion}

A well-known inverse problem from geophysics is that of full wave inversion (FWI), whereby the material properties of a medium are inferred from data on the propagation of waves in that medium \cite{mercier_designing_2025}. This problem is challenging due to ill-posedness and nonconvexity. In particular, the presence of many local minima may render standard gradient descent strategies ineffective \cite{kuldeep_full_2021}. Thus, we use the FWI problem as a test case for the LGF implementation of the ADAM algorithm, where the gradient is approximated with the surrogate. The propagation of pressure waves in a stationary medium is governed by
\begin{equation}\label{wave_continuous}
    m^2(\mathbf{x}) \pdd{u}{t} - \nabla^2 u = f(\mathbf{x},t), \quad u(\mathbf{x},0) = 0,
\end{equation}
where $u(\mathbf{x},t)$ is the pressure field, $f(\mathbf{x},t)$ is a source term, and $m(\mathbf{x})$ is the ``slowness'' of the medium, or the inverse of the wave speed. Note that the variable $t$ here is physical time rather than the pseudo-time of the optimization dynamics as in Eq. \eqref{eq:gradient_flow}. We take the computational domain to be $\Omega=[0,1]^2$ with the following boundary conditions
\begin{equation} \label{wave_bc}
\begin{aligned}
    & u(\mathbf{x}) = 0, \quad \mathbf{x} \in \partial \Omega_D, \\
    &\nabla u \cdot \mathbf{\hat n} = 0, \quad \mathbf{x} \in \partial \Omega_N,
\end{aligned}
\end{equation}

\begin{figure}[hbt!]
\centering
\includegraphics[width=1.0\textwidth]{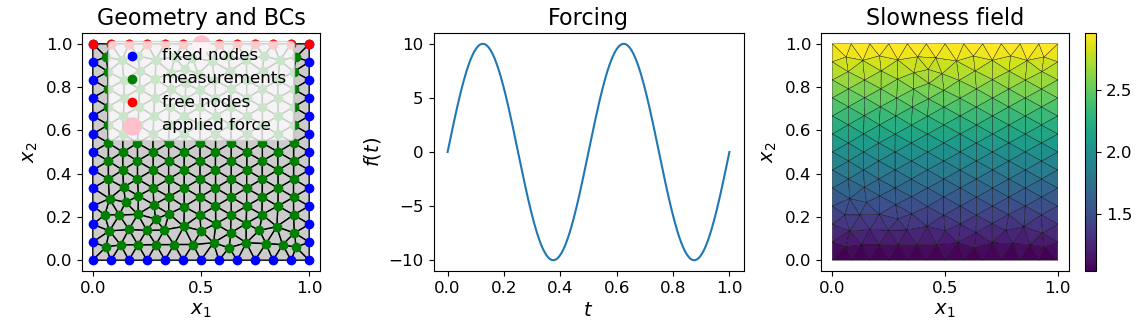}
\caption{Problem setup for the full wave inversion problem. Three of the four surfaces have homogeneous Dirichlet boundaries, while the upper horizontal surface is free. A sinusoidal point force is applied at the center of this surface. There are $n_{el}=346$ elements in the mesh, and the true slowness field is a linearly graded material.}
\label{scattering_setup}
\end{figure}

\noindent where $\partial \Omega_D$ is the Dirichlet portion of the boundary, $\partial \Omega_N$ is the Neumann portion, and $\mathbf{\hat n}$ is the outward-facing normal vector. The source term is $f(\mathbf{x},t) = \delta(x_1-1/2,x_2-1) 10 \sin(4\pi t)$. The true slowness field is that of a functionally graded material and is given by 
\begin{equation*}
    m_{\text{true}}^2(\mathbf{x}) = 1 + 2x_2.
\end{equation*}
Three of the surfaces of the geometry have zero-pressure boundary conditions, whereas the upper surface is free. See Figure \ref{scattering_setup} for the problem geometry, boundary conditions, and finite element mesh. The goal of the inverse problem is to recover the squared slowness field $m^2(\mathbf x)$ from data. Analogous to the density in the topology optimization problem, we discretize the squared slowness field as constant over elements in a finite element mesh. Thus, the optimization variables $\mathbf{a}$ represent the elemental values of the squared slowness fields. Similarly, we discretize the pressure field in space and time with
\begin{equation*}
    u(\mathbf{x},t) = \sum_{n=1}^{n_{node}} u_i(t) N_i (\mathbf{x}),
\end{equation*}
where the $N_i(\mathbf{x})$ are linear finite element functions on a triangular mesh and $n_{node}$ is the number of nodes in the mesh. We compute the Bubnov-Galerkin weak form of Eq. \eqref{wave_continuous} in space to obtain the discretized governing equation as
\begin{equation}\label{governing}
    \mathbf{M}(\mathbf{a}) \pdd{\mathbf{u}}{t} + \mathbf{K} \mathbf{u}(t) = \mathbf{F}(t), \quad \mathbf{u}(0) = \mathbf{0}.
\end{equation}

By assumption, the boundary conditions of Eq. \eqref{wave_bc} are built into the finite element discretization of the pressure field. To set up the inverse problem, we assume that the nodal displacements are observed continuously in time without noise. The measurement data is given by $\mathbf{v}(t)$ where $t \in [0,T]$. The inverse problem for the material property field is then
\begin{equation}\label{wave_opt}
\begin{aligned}
    & \underset{\mathbf{u}, \mathbf{a}}{\text{argmin }} \frac{1}{2} \int_0^T \lVert  \mathbf{u}(t) - \mathbf{v}(t) \rVert^2 dt \\
    & \text{s.t. }  \mathbf{M}(\mathbf{a}) \pdd{\mathbf{u}}{t} + \mathbf{K} \mathbf{u}(t) - \mathbf{F}(t) = \mathbf{0}, \quad \mathbf{u}(0)=\mathbf{0}.
\end{aligned}
\end{equation}

As in the previous problems, we can treat the pressure as an explicit function of the material parameters $\mathbf{a}$ through the governing equation. To compute the loss with the measurement data in Eq. \eqref{wave_opt} at a given $\mathbf{a}$, we time integrate the governing equation with central differencing, which gives rise to the update rule
\begin{equation}\label{wave_integrate}
    \mathbf{u}_{k+1} = \Delta t^2 \mathbf{M}(\mathbf{a})^{-1}\Big(  \mathbf{F}_k - \mathbf{K} \mathbf{u}_k \Big) + 2 \mathbf{u}_k - \mathbf{u}_{k-1}.
\end{equation}
We will denote a solution for pressure through Eq. \eqref{wave_integrate} as $\mathbf{u}(t;\mathbf{a})$. The constrained optimization problem of Eq. \eqref{wave_opt} can be reformulated as 
%
\begin{equation}\label{wave_obj}
\begin{aligned}
    & \underset{\mathbf{a}}{\text{argmin }} z(\mathbf{a})
    \\ & z(\mathbf{a}) = \frac{1}{2} \int_0^T \lVert \mathbf{u}(\mathbf{a},t) - \mathbf{v}(t) \rVert^2 dt .
\end{aligned}
\end{equation}

\begin{figure}[hbt!]
\centering
\includegraphics[trim = 20mm 20mm 20mm 15mm, clip, width=1\textwidth]{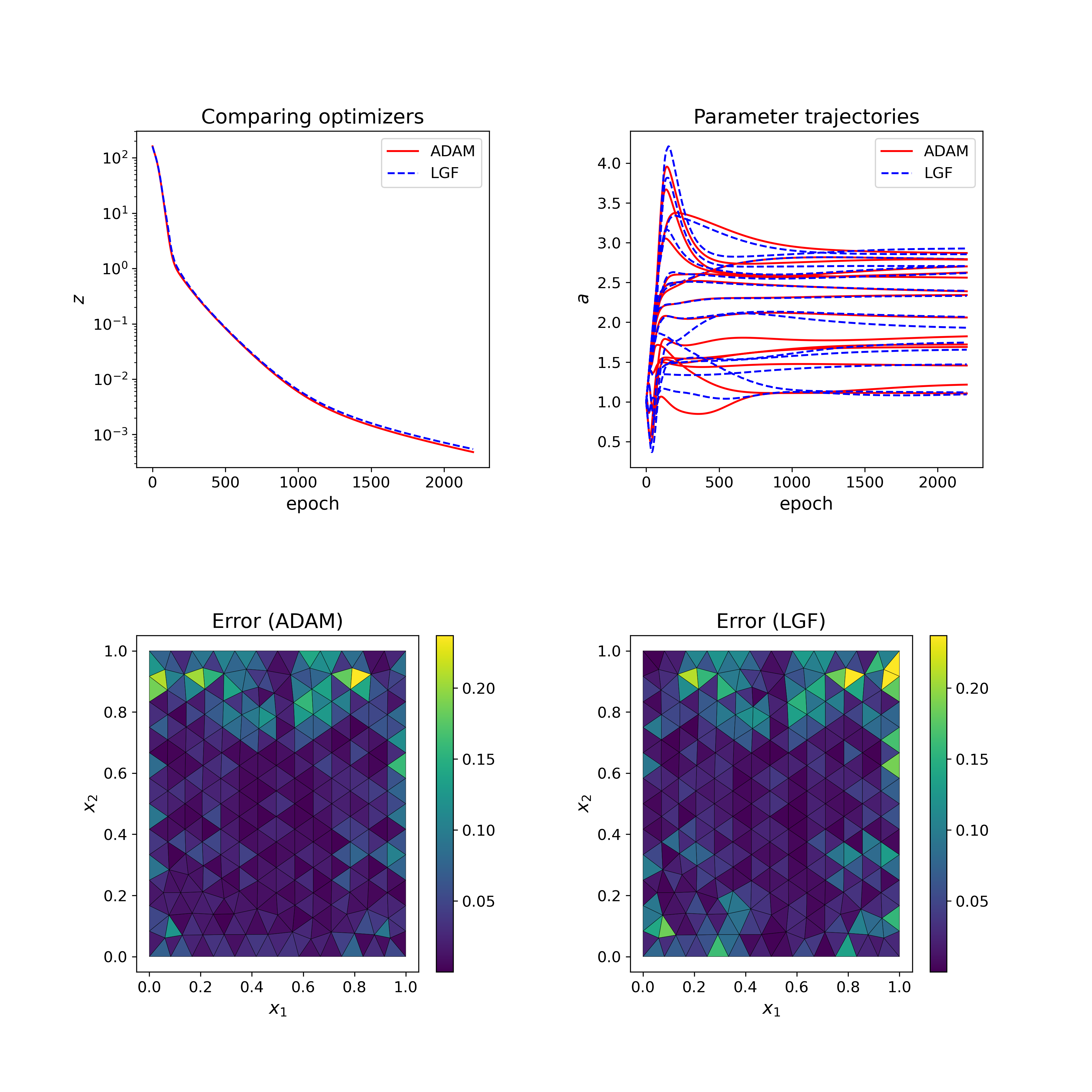}
\caption{Comparing standard ADAM optimization to the continuous time ADAM with gradient surrogate. Note that the elemental slowness parameters have very little influence on the predicted solution (and thus, the loss) as the optimizer approaches convergence. The surrogate model tracks closely with the parameter evolution from true ADAM optimization. The error reported here is the absolute value of the difference between the true and recovered slowness fields. This FWI problem is accelerated by $50\%$.}
\label{scattering_ex}
\end{figure}

In order to solve this unconstrained optimization problem, we require the gradient $\partial z / \partial \mathbf{a}$, which in turn requires $\partial \mathbf{u} / \partial \mathbf{a}$. Though adjoint state methods are a standard tool in FWI problems to obtain this gradient, we opt to use automatic differentiation in PyTorch to compute the necessary gradients of Eq. \eqref{wave_obj}. The data is generated by time integrating Eq. \eqref{governing} with the true slowness field and given source term. To avoid a complex solution field from interacting reflections, we choose the simulation time as $T=1$. We use $300$ time integration points in the central difference method of Eq. \eqref{wave_integrate}. We compare a standard ADAM optimizer to the continuous-time ADAM with the surrogate for the gradient. We take the learning rate to be $\texttt{lr}=2 \times 10^{-2}$, assume the standard Adam hyperparameters as $\beta_1 = 0.9, \beta_2 = 0.999, \epsilon=10^{-8}$, use a quadratic approximation of the gradient dynamics with a polynomial order of $\mathcal P=2$, and set $\mathcal M=30$ and $\mathcal K=20$.
Unlike the other examples, within the \texttt{STLSQ} solver, we set $\alpha = 0$ and \texttt{max\_iter~=~5} due to the high cost of this problem. We choose to run the optimization for a fixed number of epochs and verify convergence by observing the objective value, which reports the error between the measured solution and the solution predicted from the current material properties. In this example, we run the optimization for $2200$ epochs.
The pointwise error in the recovered slowness field is taken as $| a(\mathbf{x}) - m^2_{\text{true}}(\mathbf{x})|$ where $a(\mathbf{x})$ is the slowness field obtained from the chosen optimizer. See Table \ref{tab:ex4} for a summary of the problem setup and Figure \ref{scattering_ex} for the results. The final value of the error between the measurement data and predicted solution in Eq. \eqref{wave_obj} from the two optimizers is
\begin{equation*}
    z(\mathbf{a}_{\text{ADAM}}) = 4.8 \times 10^{-4}, \quad z(\mathbf{a}_{\text{LGF}}) = 5.5 \times 10^{-4}.
\end{equation*}

We remark that, near convergence, the elemental slownesses have a very small influence on the predicted spacetime solution field. This means that the corresponding objective gradient terms become small and the slowness field converges slowly. Accordingly, we do not necessarily expect perfect recovery of the material parameters. However, we do demand that the ADAM variant of the LGF optimizer has comparable performance to the standard implementation of ADAM. The relative mean-squared errors with the exact solution from the two optimizers are given by 
\begin{equation*}
    \frac{ \int_{\Omega} | a_{\text{ADAM}}(\mathbf{x}) - m^2_{\text{true}}(\mathbf{x})|^2 d\Omega} {\int_{\Omega} | m^2_{\text{true}}(\mathbf{x}) |^2 }= 7 \times 10^{-4}, \quad \frac{ \int_{\Omega} | a_{\text{LGF}}(\mathbf{x}) - m^2_{\text{true}}(\mathbf{x})|^2 d\Omega }{\int_{\Omega} | m^2_{\text{true}}(\mathbf{x}) |^2}=9 \times 10^{-4}.
\end{equation*}
In the context of a dynamic inverse problem, each time the surrogate is used in place of evaluating the gradient, we avoid backpropagating through the time integration, which finds the solution field at the current material properties. Both the surrogate and ADAM achieve comparable performance on the FWI problem, driving the objective down $5$ orders of magnitude and leading to accurate recovery of the true slowness field.

\begin{table}
    \centering
    \begin{tabular}{|c|c|}
        \hline
        \textbf{Problem parameter} &  \textbf{Setting} \\
        \hline
        $\#$ optimization variables & $346$  \\
        $\texttt{lr} $& $2 \times 10^{-2}$  \\
        Retrain interval, $\mathcal M$ & $30$ \\
        History size, $\mathcal K$ & $20$\\
        Polynomial order, $\mathcal P$ & $2$ \\
        $\texttt{epochs}$ & $2200$ \\
        Truncation rank, $r$ & N/A\\
         \hline
    \end{tabular}
    \caption{Problem parameters for full wave inversion problem.}
    \label{tab:ex4}
\end{table}

\begin{figure}[hbt!]
\centering
\includegraphics[trim = 20mm 10mm 10mm 15mm, clip, width=1\textwidth]{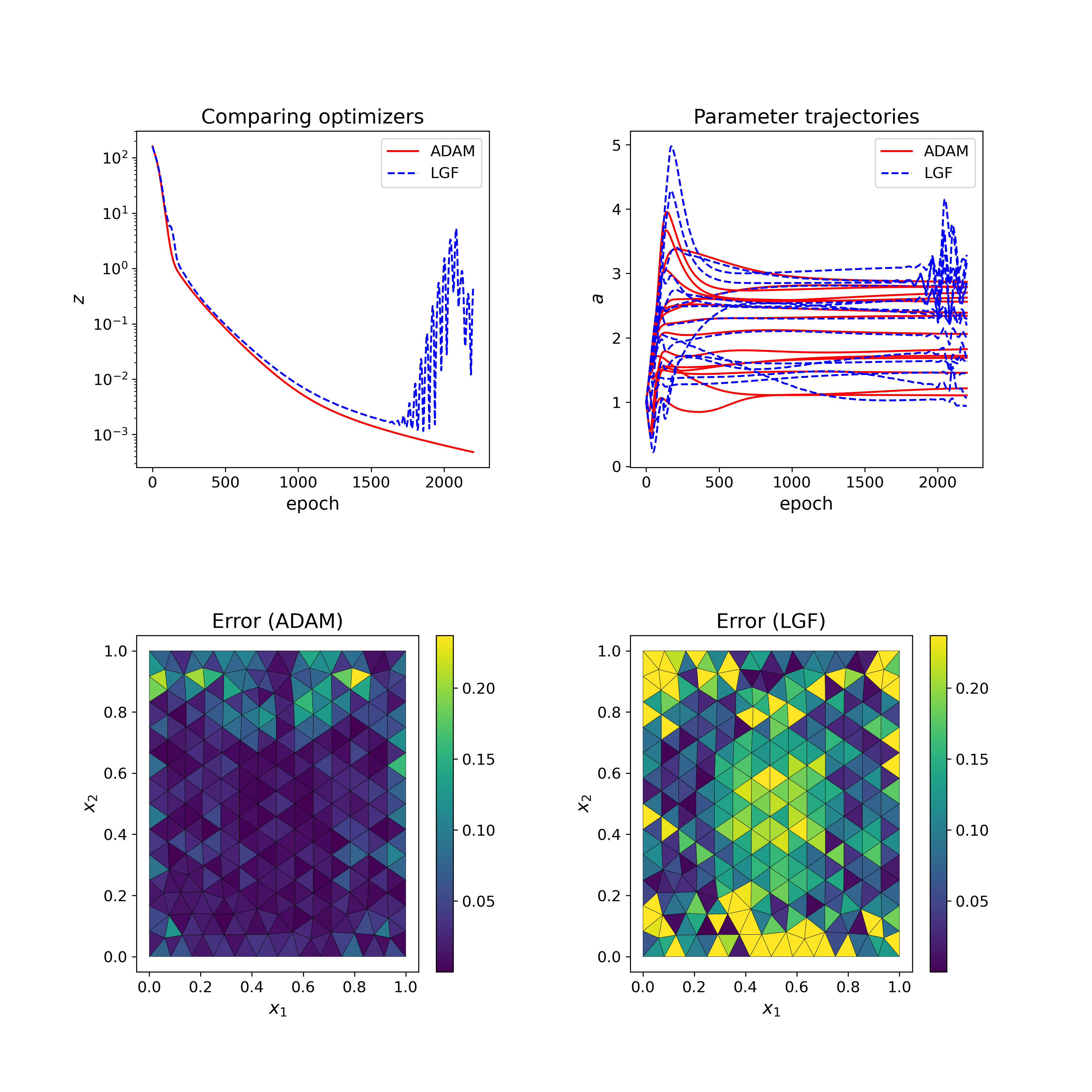}
\caption{Comparing standard ADAM optimization to the continuous-time ADAM with gradient surrogate and $100\%$ acceleration. Towards the end of the optimization, the surrogate consistently increases the objective function, indicating that the optimization dynamics are inaccurately approximated. The differences in the recovered slowness fields between the two methods are now pronounced.}
\label{scattering_ex2}
\end{figure}

As in the previous example, we experiment with more aggressive acceleration of the FWI problem. With all other problem parameters equivalent, we set the retrain interval at $\mathcal M=40$, corresponding to $100 \%$ acceleration. The relative mean-squared error obtained by the LGF is $2.3 \times 10^{-3}$. See Figure \ref{scattering_ex2} for the results. In the first $1750$ epochs of optimization, the surrogate drives down the loss monotonically, but then fails to capture the optimization dynamics closer to convergence. In these cases, the errors incurred by the surrogate need to be corrected by ADAM, which indicates that the surrogate is not providing any benefit. Again, this example illustrates the importance of appropriately tuning the history size and retrain interval hyperparameters.

\subsection{Deep Ritz method with ADAM optimization}

\begin{figure}[hbt!]
\centering
\includegraphics[width=1\textwidth]{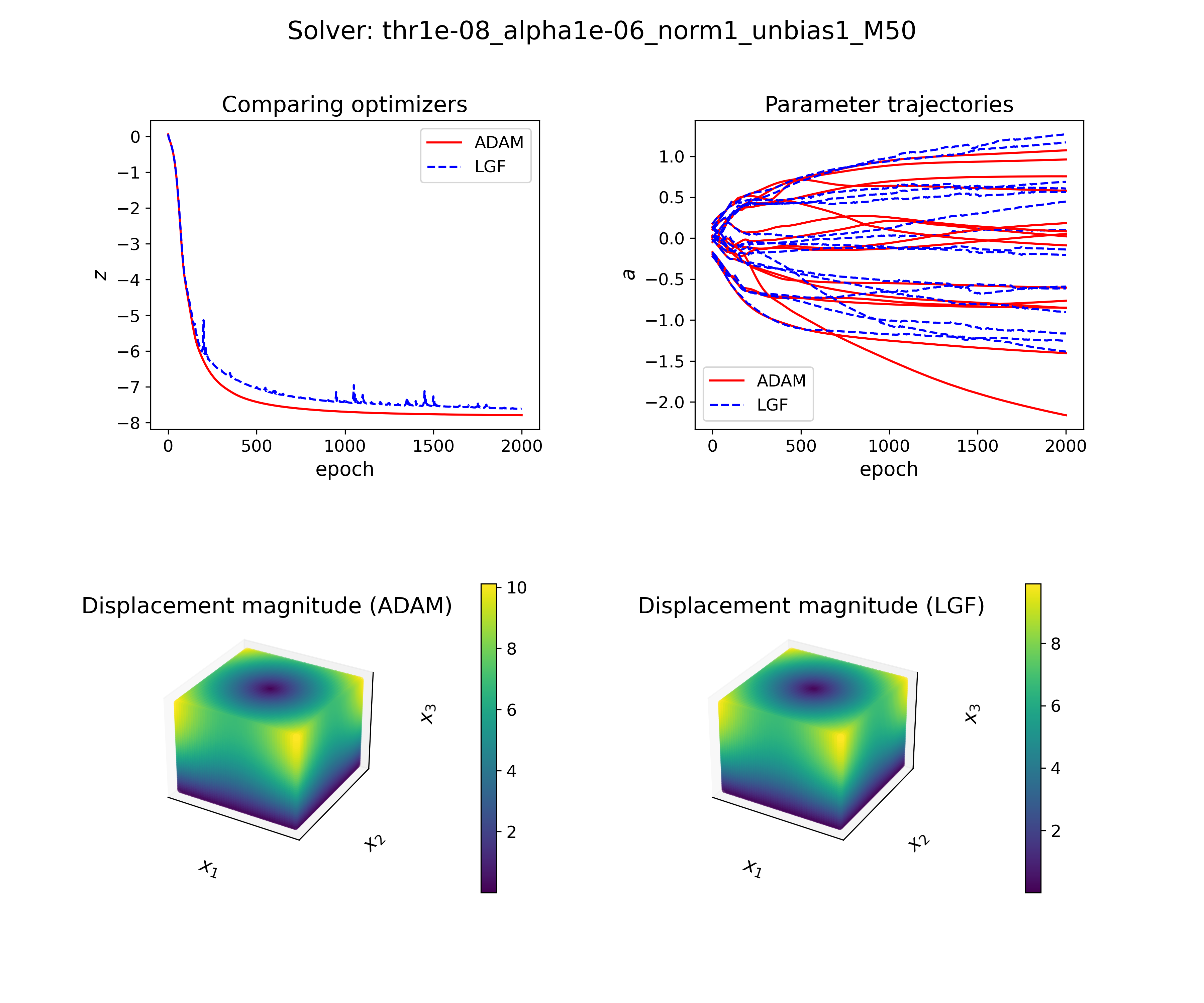}
\caption{Comparing standard ADAM optimization with the LGF version of ADAM. We save $15$ out of every $50$ evaluations of the gradient of the loss function using the surrogate model, corresponding to an acceleration of $43\%$.}
\label{elasticity}
\end{figure}

As a final test of our method, we show an example involving neural network training with ADAM optimization. In particular, we train a neural network to represent the displacement field of a three-dimensional linearly elastic structure. In Eq. \eqref{to_continuous}, we work with the strong form of the governing equation for linear elasticity. An alternative formulation of static equilibrium is that of the principle of minimum total potential energy, which says that the total potential energy functional is minimized. In the continuous setting, this reads
\begin{equation}\label{total_potential}
    \Pi = \int_{\Omega} \frac{1}{2} \sigma_{ij} (\mathbf{x})\epsilon_{ij}(\mathbf{x}) - b_i(\mathbf{x}) u_i(\mathbf{x}) d\Omega - \int_{\partial \Omega} t_i(\mathbf{x}) u_i (\mathbf{x}) dS, \quad \delta \Pi = 0,
\end{equation}
where $t_i$ is the surface traction and $\Pi$ is the total potential energy. The Deep Ritz Method \cite{e_deep_2017} approximates a solution to this problem by discretizing the displacement field with a neural network and minimizing the energy functional in terms of the neural network parameters. Because Eq. \eqref{total_potential} assumes that the Dirichlet boundary conditions are satisfied by the displacement automatically, we choose to build in boundary conditions to the neural network discretization with a distance function method \cite{sukumar_exact_2022}. The displacement field is thus discretized as 
\begin{equation}\label{distance}
    \mathbf{u}(\mathbf{x}; \mathbf a) = D(\mathbf{x}) \mathcal{ N}(\mathbf{x};\mathbf a),
\end{equation}
\begin{table}
    \centering
    \begin{tabular}{|c|c|}
        \hline
        \textbf{Problem parameter} &  \textbf{Setting} \\
        \hline
        $\#$ optimization variables & $825$  \\
        $\#$ integration points & $125000$  \\
        $\texttt{lr} $& $2.5 \times 10^{-3}$  \\
        Retrain interval, $\mathcal M$ & $50$ \\
        History size, $\mathcal K$ & $35$\\
        Polynomial order, $\mathcal P$ & $1$ \\
        $\texttt{epochs}$ & $2000$ \\
        Truncation rank, $r$ & N/A\\
         \hline
    \end{tabular}
    \caption{Problem parameters for Deep Ritz solution to linearly elastic deformation.}
    \label{tab:ex5}
\end{table}
where $\mathcal{ N}: \mathbb{R}^3 \rightarrow \mathbb{R}^3$ is a multi-layer perceptron neural network, $D(\mathbf{x})$ is a distance function which is zero along the Dirichlet portion of the boundary, and $\mathbf a$ are the optimization variables, which are the parameters of the neural network. Note that the form of Eq. \eqref{distance} assumes homogeneous Dirichlet boundary conditions. In our example, we take the computational domain to be the unit cube $\Omega=[0,1]^3$, the $x_3=0$ surface to be clamped, and the remaining portion of the boundary to be traction-free. The distance function is taken as $D(\mathbf{x})=x_3$, and we use a two-hidden-layer neural network with $25$ hidden units per layer and hyperbolic tangent activation functions. The integration grid used to form the total potential energy is $50 \times 50 \times 50$ evenly spaced points. The volumetric force is taken to be that of a torsional displacement field, i.e.,
\[
\mathbf{b}(\mathbf{x})
= b_0
\begin{bmatrix}
-(x_2 - \tfrac{1}{2}) & (x_1 - \tfrac{1}{2}) & 0
\end{bmatrix}^{T}.
\]
where $b_0=10$ controls the magnitude of the forcing. We compare ADAM optimization against the LGF optimizer using the ADAM setting with a linear approximation of the gradient, history size $\mathcal K=35$, and retrain interval $\mathcal M=50$. A summary of the parameter settings is given in Table \ref{tab:ex5}. See Figure \ref{elasticity} for the results. The final energy values after $2000$ epochs of true ADAM optimization and ADAM optimization with a linear surrogate for the gradient are given by
\begin{equation*}
   \Pi_{\text{ADAM}} = -7.79, \quad \Pi_{\text{LGF}} = -7.60. 
\end{equation*}

This corresponds to a normalized difference in predicted displacement fields given by 
\begin{equation*}
    \frac{ \int_{\Omega} \lVert \mathbf{u}_{\text{ADAM}}(\mathbf{x}) - \mathbf{u}_{\text{LGF}}(\mathbf{x})\rVert^2 d\Omega}{ \int_{\Omega}  \lVert \mathbf{u}_{\text{ADAM}}(\mathbf{x})\rVert^2 d\Omega} = 1.3 \times 10^{-3}.
\end{equation*}

Note that the optimization variable trajectories from the two optimizers do not coincide as neatly as the cases of the previous examples. However, this does not indicate decreased performance of the surrogate model, as seen in close agreement of the final values of the loss and the predicted displacement fields. Non-uniqueness in the dependence of the neural network predictions on the optimization variables suggests that, in this case, the loss trajectory is a better indicator of performance than the agreement of the parameter trajectories from the two optimizers. The true ADAM dynamics lead to a lower converged energy value, but this leads only to slight differences in the displacement fields.

\begin{figure}[hbt!]
\centering
\includegraphics[width=1\textwidth]{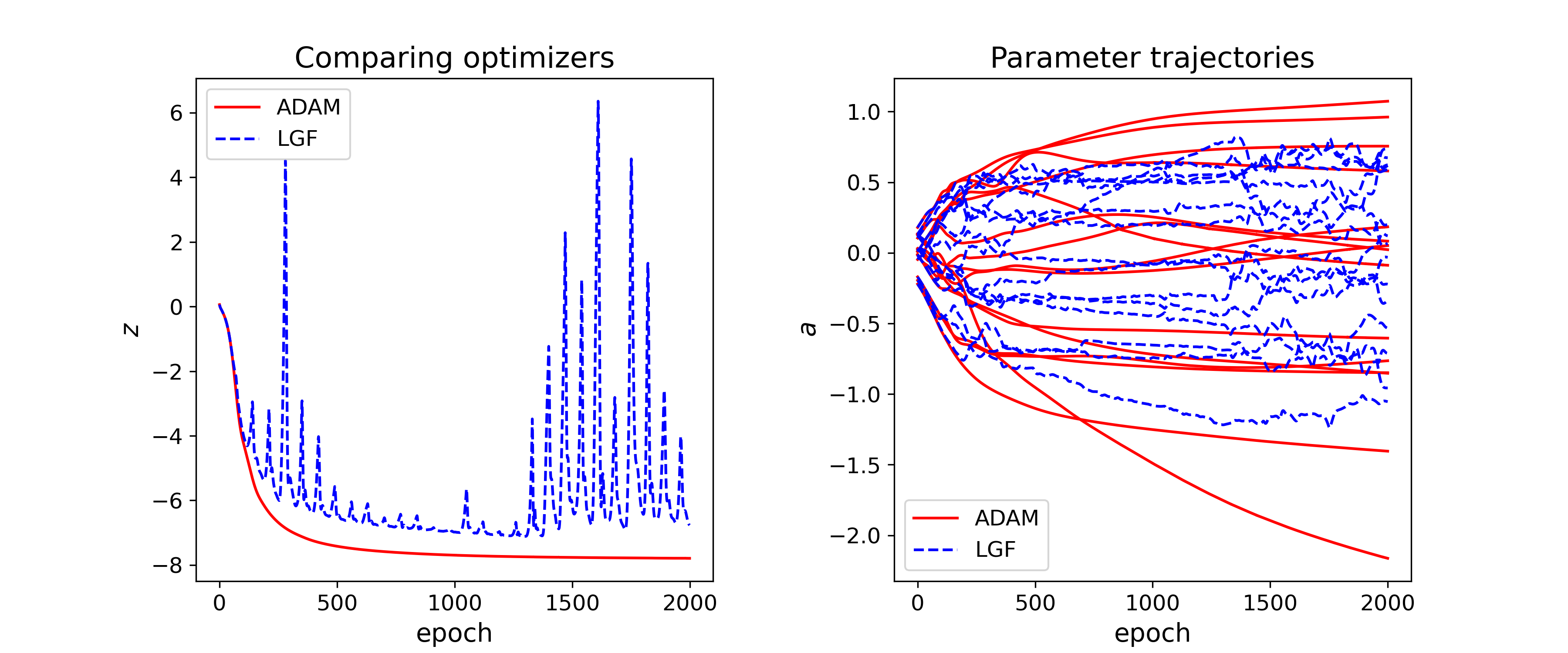}
\caption{When the acceleration is increased to $100\%$, the surrogate model does not accurately reflect the loss landscape. When switched on, it consistently drives the total potential energy up, providing no benefit to the optimization process.}
\label{elasticity3}
\end{figure}

With this Deep Ritz example, we again gain insight into the pitfalls of choosing too large an acceleration value. Keeping all problem parameters constant, we increase the acceleration to $100\%$ by setting $\mathcal M=70$. See Figure \ref{elasticity3} for the results. The surrogate model consistently drives up the loss, and these errors need to be reversed by the true ADAM dynamics. We see that the surrogate model provides no benefit when the acceleration is too large.

\begin{figure}[hbt!]
\centering
\includegraphics[width=1\textwidth]{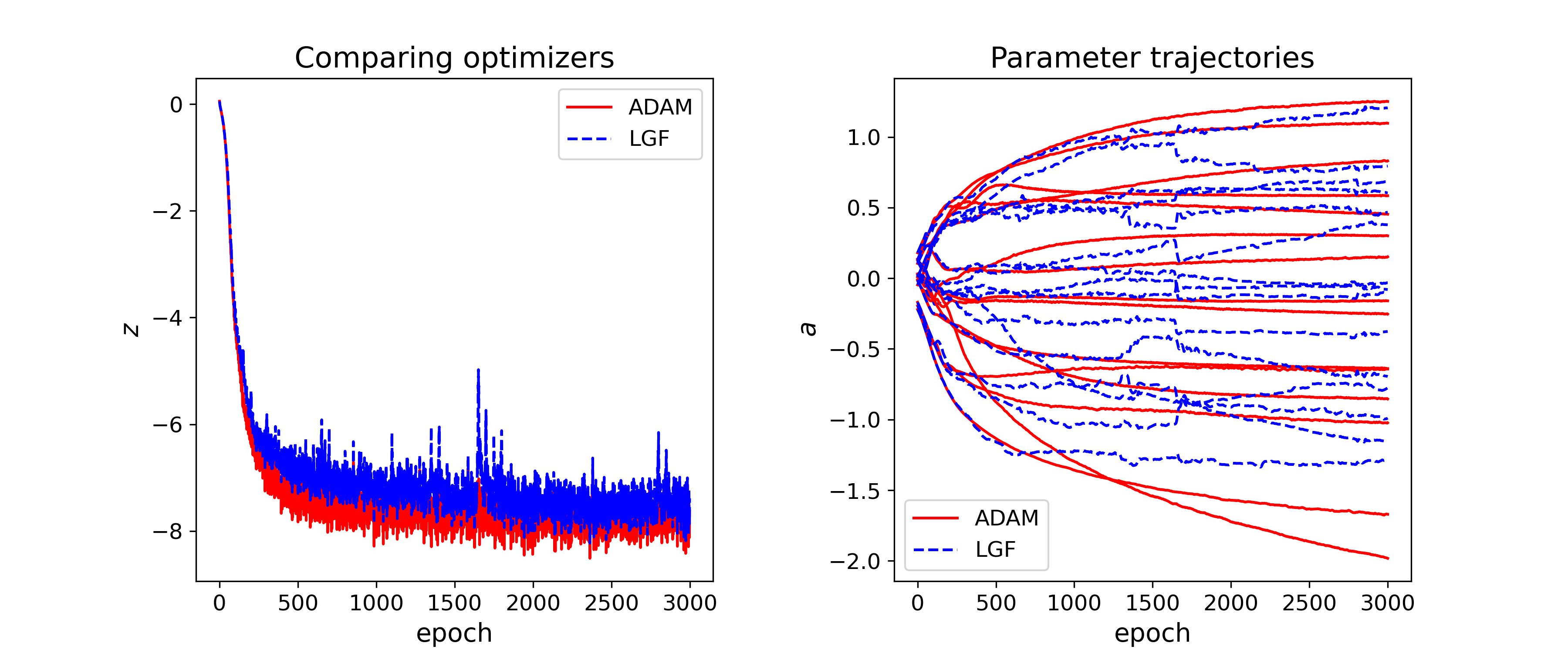}
\caption{Comparing the two optimizers with a Monte Carlo approximation of the total potential energy. The loss convergence indicates noisy gradients, yet the surrogate model still captures relevant features of the optimization dynamics. Once again, we save $15$ out of every $50$ evaluations of the gradient of the loss function using the surrogate model, corresponding to an acceleration of $43\%$.}
\label{elasticity_stochastic}
\end{figure}

We also test our method on a stochastic version of the strain energy objective given in Eq. \eqref{total_potential}. Such a test demonstrates that the learned gradient flow can perform on noisy gradients from machine learning problems where the training data is mini-batched, as is often the case in practice. The analogue of mini-batching in the physics-informed setting is to use Monte Carlo integration to form the total potential energy at each optimization epoch. Namely, we approximate the energy objective as 
\begin{equation}\label{stochastic_potential}
    \hat \Pi \approx \frac{1}{B} \sum_{i=1}^B \frac{1}{2} \sigma_{ij}(\mathbf{x}_i) \epsilon_{ij}(\mathbf{x}_i) - b_i(\mathbf{x}_i) \hat u_i(\mathbf{x}_i), \quad \mathbf{x}_i \overset{\text{i.i.d.}}{\sim} \mathcal{U}(0,1),
\end{equation}
where $B$ is the size of the mini-batched integration. We take the integration batch to be $B=5000$, which is $4\%$ of the size of the full-batch integration grid, and run the optimization for $3000$ epochs with all other problem and optimization variables set at the same values. We remark that the stochastic integration routine requires $1.5 \times 10^7$ evaluations of the neural network for the displacement, whereas the full-batch integration requires $2.5 \times 10^8$ evaluations. Mini-batching integration can be an effective strategy to accelerate neural network-based physics problems, as noted in \cite{rowan_variational_2025}. See Figure \ref{elasticity_stochastic} for the results of training with the stochastic objective. The final values of the energy for the two optimizers are
\begin{equation*}
   \hat \Pi_{\text{ADAM}} = -7.62, \quad \hat \Pi_{\text{LGF}} = -7.39. 
\end{equation*}

The normalized difference in the displacement fields obtained by minimizing the stochastic objective with the two optimizers is 
\begin{equation*}
    \frac{ \int_{\Omega} \lVert \mathbf{\hat u}_{\text{ADAM}}(\mathbf{x}) - \mathbf{\hat u}_{\text{LGF}}(\mathbf{x})\rVert^2 d\Omega}{ \int_{\Omega}  \lVert \mathbf{\hat u}_{\text{ADAM}}(\mathbf{x})\rVert^2 d\Omega} = 3.1 \times 10^{-3}.
\end{equation*}

Again, we note that the trajectories of the two sets of neural network parameters do not overlap as closely in the Deep Ritz problem as in previous examples. This is a consequence of the non-uniqueness of neural network approximations. 

\begin{table}[htb]
    \centering
    \begin{tabular}{|c|c|c|}
        \hline
        \textbf{Optimization problem} &  \textbf{Method} & \textbf{Acceleration} \\
        \hline
        Dynamic inverse problem & Gradient descent & $200 \% $ \\
        Topology optimization & Gradient descent &  $250 \% $ \\
        Nonlinear heat conduction & Newton's method & $33 \% $ \\
        Full wave inversion & ADAM & $50 \%$ \\
        Deep Ritz training & ADAM (full-batch) & $43 \%$ \\ 
        Deep Ritz training & ADAM (mini-batch) & $43 \%$ \\ 
        \hline
    \end{tabular}
    \caption{Summary of the optimization problems studied, the optimization method, and the acceleration observed from using the data-driven surrogate. The acceleration is computed as $100 \times (\mathcal{M} / \mathcal{K} - 1)$.}
    \label{tab:acceleration_summary}
\end{table}

\section{Conclusion}

In this work, we connected two key ideas: gradient flows and data-driven equation discovery. During training, optimization variables evolve according to an underlying gradient flow, which can be viewed as a dynamical system.
Separately, equation discovery methods seek to approximate the behavior of dynamical systems from their observations as a means to generate future predictions.
We combined these perspectives by learning a dynamical system that approximates the gradient flow of the optimizer.
This learned gradient flow then predicts the trajectory of the optimization variables, resulting in computational savings by avoiding evaluations of the true gradient dynamics, which often involve expensive PDE solves.
We initially derived this method for gradient descent, but generalized it to Newton's method and ADAM.
We experimentally demonstrated the utility of this idea across a range of mechanics problems, including topology optimization, a nonlinear heat conduction inverse problem, full wave inversion, and the Deep Ritz method.
Across these examples, our learned gradient flow reduced the number of gradient evaluations by factors ranging from $33 \%$ up to $250 \%$, as summarized in Table \ref{tab:acceleration_summary}.

While these results demonstrate clear benefits, several challenges and open questions remain.
A particularly important direction is developing principled guidance for selecting the algorithmic parameters that govern when and how the surrogate is learned and deployed.
While we provided cursory examples of changing the retrain interval $\mathcal{M}$ at the end of Sec. \ref{sec:ex:heat-newton} and \ref{sec:ex:wave-inversion}, the selection of this algorithmic parameter warrants further study well beyond what we have provided here.
Rather than alternating between the surrogate and true training dynamics on a fixed schedule, a possible alternative is to design explicit criteria for determining when retraining is necessary, thus automating the choice of acceleration.
Similarly, some methods may automate selection of the model form ($\mathcal{P}$ and $r$ in our case) based on a train-validation split.
Many advances in dynamical systems discovery may also be applied to this scenario, opening a multitude of possible improvements.
In particular, one could modify this approach with an online version of SINDy that operates on streaming data, or perhaps decouple the model-building step from the training so it can run in parallel or on a separate device.
Similarly, the (optional) dimensionality reduction step could be replaced with an online equivalent.
Another direction could investigate incorporating components of other optimizers into this paradigm, such as the Hessian approximation afforded by LBFGS or various improvements to ADAM \cite{reddi_convergence_2019, loshchilov_decoupled_2019}.
The LGF paradigm is also broadly applicable to iterative optimization methods in general, including gradient-free methods, constituting another possible research direction.
On the practical side, several considerations merit further attention: developing implementations for distributed optimization problems, automating hyperparameter selections, and reducing memory overhead.

Finally, it is important to recognize that this approach is, unsurprisingly, not universally applicable.
While some of the aforementioned techniques may help to generalize to high-dimensional optimization problems, there will always be problems whose training dynamics cannot be effectively captured by a surrogate learned from limited, original optimization iterations.
In particular, LGF may struggle in constrained optimization settings where the algorithm intermittently switches between different constraint regimes (e.g., via projection steps, feasibility restoration, active-set updates, or complementarity conditions).
Such switching can introduce non-smooth behavior in the update direction, producing abrupt changes in the observed trajectories that an autonomous, smooth surrogate learned from a short window of iterations is unlikely to extrapolate reliably.


\section*{Acknowledgment}
 This work was supported by the Department of Energy's National Nuclear Security Administration under Award Numbers DE-NA0003962 and DE-NA0003968. C. Rowan's work was funded by the National Defense Science and Engineering Graduate Fellowship (NDSEG) through the Department of Defense (DOD) and the Army Research Office (ARO).



\clearpage
\appendix
\renewcommand{\thesection}{\Alph{section}} 
\titleformat{\section}
  {\normalfont\Large\bfseries}
  {Appendix \thesection}{1em}{}
  

\section{ADAM ODE Verification}
\label{appendix:adam_ode_verification}

Given the ODE system for ADAM in Eq. \eqref{eq:adam_ode_with_delta_t}, we show how we can recover the standard ADAM update rules in Eq. \eqref{eq:discrete_adam:param_update} -- \eqref{eq:discrete_adam:unbiased_update} through time discretization.
We repeat the continuous system below for convenience:
\begin{equation}
    \begin{aligned}
        \frac{\dd \ab}{\dd t} (t) &= - \frac{\frac{\mb (t)}{1 - \beta_{1}^{t/\eta}}}{\sqrt{\frac{\vbb (t)}{1 - \beta_{2}^{t/\eta}}} + \epsilon} , \\
        \frac{\dd \mb}{\dd t} (t) &=  \frac{1}{\eta} (1 - \beta_1) \left (  \pd{}{\mathbf{a}} z (\ab_{k})  - \mb(t) \right ), \\
        \frac{\dd \vbb}{\dd t} (t) &= \frac{1}{\eta} (1-\beta_2)  \left ( \left ( \pd{}{\mathbf{a}} z (\ab_{k}) \right ) ^{2} - \vbb (t) \right ) .
    \end{aligned}
    \tag{\ref{eq:adam_ode_with_delta_t}} \label{eq:adam_ode_with_delta_t_reprint}
\end{equation}
We will use a semi-implicit method to discretize this system, applying forward Euler for $\mb$ and $\vbb$ and backward Euler for $\ab$.
We apply the forward Euler update at $t_{k} = \eta \cdot k$ with a step size of $\eta$, denoting $\mb_{k} = \mb (t_k)$, $\vbb_{k} = \vbb (t_k)$, and $\ab_{k} = \ab(t_k)$, this gives
\begin{align*}
    \frac{\mb_{k+1} - \mb_{k}}{\eta} &=  \frac{1}{\eta} (1 - \beta_1) \left ( \pd{}{\mathbf{a}} z (\ab_{k}) - \mb_{k} \right ), \\
    \frac{\vbb_{k+1} - \vbb_{k}}{\eta} &= \frac{1}{\eta} (1-\beta_2) \left ( \left ( \pd{}{\mathbf{a}} z (\ab_{k}) \right ) ^{2} - \vbb_{k} \right ) .
\end{align*}
Expanding the $(1 - \beta_1)$ and $(1 - \beta_2)$ terms with $\mb_k$ and $\vbb_k$ and simplifying gives the update rule for the biased moments from Sec. \ref{sec:learned_adam}:
\begin{equation}
    \begin{aligned}
        \mb_{k+1} &= \beta_1 \mb_{k} + (1 - \beta_1) \, \pd{}{\mathbf{a}} z (\ab_{k}) , \\
        \vbb_{k+1} &= \beta_2 \vbb_{k} + (1 - \beta_2) \left (  \pd{}{\mathbf{a}} z (\ab_{k}) \right )^2 .
    \end{aligned}
    \tag{\ref{eq:discrete_adam:unbiased_update}} \label{eq:discrete_adam:unbiased_update_reprint}
\end{equation}
 
Next, returning to the ODE system in Eq. \eqref{eq:adam_ode_with_delta_t_reprint}, applying backward Euler for the update to $\ab$, and noting $t_{k+1} = (k+1) \cdot \eta$,
\begin{equation}
    \frac{\ab_{k+1} - \ab_{k}}{\eta} =
    - \frac{
        \frac{
            \mb_{k+1}
        }
        {
            1 - \beta_{1}^{k+1}
        }
    }
    {
        \sqrt{
            \frac{
                \vbb_{k+1}
            }
            {
                1 - \beta_{2}^{k+1}
            }
        } + \epsilon
    } . 
    \label{eq:appendix:adam_derivation:nested}
\end{equation}
For the next step, recall the unbiased moment formulas:
\begin{equation}
    \hat{\mb}_{k+1} = \frac{\mb_{k+1}}{1 - \beta_{1}^{k+1}}, \qquad 
    \hat{\vbb}_{k+1} = \frac{\vbb_{k+1}}{1 - \beta_{2}^{k+1}}.
    \tag{\ref{eq:discrete_adam:biased_moments}} \label{eq:discrete_adam:biased_moments_reprint}
\end{equation}  

The nested fractions in Eq. \eqref{eq:appendix:adam_derivation:nested} correspond to the unbiased estimators of Eq. \eqref{eq:discrete_adam:biased_moments_reprint}, and through that substitution, we have recovered the original update rule in Eq. \eqref{eq:discrete_adam:param_update}:
$$
    \ab_{k+1} = \ab_{k} 
    - \eta \cdot \frac{
        \hat{\mb}_{k+1}
    }
    {
        \sqrt{\hat{\vbb}_{k+1}} + \epsilon
    } . \\
$$
Thus, through a semi-implicit discretization of the ODE system in Eq. \eqref{eq:adam_ode_with_delta_t_reprint}, we recover the discrete ADAM update rule.

\end{document}